\documentclass[aos,preprint]{imsart}
\newcommand*{\Scale}[2][4]{\scalebox{#1}{$#2$}}%
\setattribute{journal}{name}{}
\usepackage[hyphens]{url}
\RequirePackage[OT1]{fontenc}
\RequirePackage{amsthm,amsmath, amsfonts, amssymb, graphics, graphicx, subfigure, color, multirow, amscd}
\usepackage{datetime}
\RequirePackage{natbib}
\RequirePackage[colorlinks,citecolor=blue,urlcolor=blue,breaklinks]{hyperref}
\usepackage{graphicx}
\usepackage{algpseudocode}
\usepackage{algorithm}% http://ctan.org/pkg/algorithms
\usepackage{comment}

\newtheorem{conj}{Conjecture}

\newtheorem{lem}{Lemma}
\newtheorem{cor}{Corollary}

\newtheorem{thm}{Theorem}
\newtheorem{prop}{Proposition}

\newtheorem{rmk}{Remark}
\newtheorem{definition}{Definition}

\usepackage{appendix}

\newcommand{\vxbarbar}{ \overline{\overline{\vx}} }

\newcommand{\ttE}{\mathtt{E}}

\newcommand{\bbE}{\mathbb{E}}
\newcommand{\bbP}{\mathbb{P}}
\newcommand{\bbR}{\mathbb{R}}
\newcommand{\bbS}{\mathbb{S}}

\newcommand{\vm}{\boldsymbol{m}}

\newcommand{\vu}{\boldsymbol{u}}

\newcommand{\vw}{\boldsymbol{w}}
\newcommand{\vx}{\boldsymbol{x}}

\newcommand{\vz}{\boldsymbol{z}}

\newcommand{\vA}{\boldsymbol{A}}
\newcommand{\vB}{\boldsymbol{B}}

\newcommand{\vD}{\boldsymbol{D}}
\newcommand{\vE}{\boldsymbol{E}}
\newcommand{\vF}{\boldsymbol{F}}
\newcommand{\vG}{\boldsymbol{G}}

\newcommand{\vI  }{\mathbf{I}}

\newcommand{\vK}{\boldsymbol{K}}

\newcommand{\vQ}{\boldsymbol{Q}}

\newcommand{\vV}{\boldsymbol{V}}
\newcommand{\vW}{\boldsymbol{W}}

\newcommand{\vZ}{\boldsymbol{Z}}

\newcommand{\bbeta}{\boldsymbol{\beta}}

\newcommand{\bgamma}{\boldsymbol{\gamma}}

\newcommand{\bLambda}{\boldsymbol{\Lambda}}

\newcommand{\bSigma}{\boldsymbol{\Sigma}}

\newcommand{\bepsilon}{\boldsymbol{\epsilon}}

\newcommand\independent{\protect\mathpalette{\protect\independenT}{\perp}}
\def\independenT#1#2{\mathrel{\rlap{$#1#2$}\mkern2mu{#1#2}}}
\newcommand{\epf}{\hfill $\Box$}

\begin{document}
\begin{frontmatter}

\title{On  the optimality of sliced inverse regression in high dimensions }
\runtitle{SIR-minimax }
\thankstext{t1}{Lin's research is supported by the Center of Mathematical Sciences and Applications at Harvard University. Liu's research is supported by the NSF Grant DMS-1120368 and NIH  Grant R01 GM113242-01 }

\begin{aug}
  	\author{\fnms{Qian} \snm{Lin}\thanksref{m1}\ead[label=e1]{qianlin@cmsa.fas.harvard.edu}\thanksref{t1}}
  	\author{\fnms{Xinran} \snm{Li}\thanksref{m1}\ead[label=e2]{xinranli@fas.harvard.edu}\thanksref{t1}}  
  	\author{\fnms{Dongming} \snm{Huang}\thanksref{m1}\ead[label=e3]{dhuang01@g.harvard.edu}\thanksref{t1}} 
\and    
	\author{\fnms{Jun S. }\snm{Liu}\thanksref{m1}\ead[label=e4]{jliu@stat.fas.harvard.edu}\thanksref{t1}}

   \runauthor{Q. Lin, X. Li, D. Huang and J. S. Liu}
  \affiliation{Harvard University\thanksmark{m1} }
\address{Qian Lin\\
    Center of Mathematical Sciences\\
    and Applications \\
    Harvard University\\
    20 Garden Street \\
    Cambridge, MA 02138  \\
    USA\\
    \printead{e1}}
\address{Xinran Li\\
	Harvard Statistics Department\\
	Science Center 7th floor\\
	One Oxford Street     \\
	Cambridge, MA 02138-2901   \\
    USA\\
    \printead{e2}}       
\address{Xinran Li\\
	Harvard Statistics Department\\
	Science Center 7th floor\\
	One Oxford Street     \\
	Cambridge, MA 02138-2901  \\ 
    \printead{e3}}       
\address{Jun S. Liu\\
	Harvard Statistics Department\\
	Science Center 7th floor\\
	One Oxford Street     \\
	Cambridge, MA 02138-2901\\   
    USA\\
    \printead{e4}}        
\end{aug}

\date{\today}

\begin{abstract}
The central subspace of a pair of random variables $(y,\vx)  \in \mathbb{R}^{p+1}$ is the minimal subspace $\mathcal{S}$ such that $y\independent \vx | P_{\mathcal{S}}\vx$. In this paper, we consider the minimax rate of estimating the central space of the multiple index models $y=f(\bbeta_{1}^{\tau}\vx,\bbeta_{2}^{\tau}\vx,...,\bbeta_{d}^{\tau}\vx,\epsilon)$ with at most $s$ active predictors where $\vx \sim N(0,\vI_{p})$. 
%Sliced inverse regression (SIR) is arguably one of the most commonly used algorithms for estimating the central space. 
We first introduce a large class of models depending on the smallest non-zero eigenvalue $\lambda$ of $var(\bbE[\vx|y])$, over which we show that an aggregated estimator based on the SIR procedure converges at rate $d\wedge((sd+s\log(ep/s))/(n\lambda))$.
We then show that this rate is optimal in two scenarios: the single index models; and the multiple index models with fixed central dimension $d$ and fixed $\lambda$.  By assuming a technical conjecture, we can show that this  rate is also optimal for multiple index models with bounded dimension of the central space. 
We believe that these (conditional) optimal rate results bring us meaningful insights of general SDR problems in high dimensions.

\iffalse
$(sd+s\log(ep/s))/(n\lambda)$

For single index models, we show under mild conditions that the optimal rate $s\log(ep/s)/(n\lambda)$ can be achieved by the SIR procedure with the help of an aggregated estimator. 
A much simpler  and computational efficient algorithm, DT-SIR, can achieve this optimal rate if the sparsity $s=O(p^{1-\delta})$ for some positive $\delta>0$. 
When the dimension $d$ of the central subspace is bounded, by assuming a technical conjecture, we can also prove that the aggregated SIR estimator achieves the optimal convergence rate $(sd+s\log(ep/s))/(n\lambda)$. Finally, the convergence rate $(sd+s\log(ep/s))/(n\lambda)$ also holds for divergent $d$ and we expect it to be optimal. 
\fi

\end{abstract}

\begin{keyword}%[class=AMS]
\kwd{sufficient dimension reduction}
\kwd{optimal rates}
\kwd{sliced inverse regression}
\kwd{semi-definite positive programming}
\end{keyword}
\end{frontmatter}

%\input{introduction}

%\input{mainresults}
%\today,
%\currenttime
\section{Introduction}
Because of rapid advances of information technologies in recent years, it has become a  common problem for data analysts that the dimension ($p$) of data is much larger than the sample size ($n$), i.e., the {\it `large $p$, small $n$ problem'}.
For these problems, variable selection and dimension reductions are often indispensable first steps. In early 1990s, a fascinating supervised dimension reduction method, the sliced inverse regression (SIR) \citep{li1991sliced},  was proposed to model univariate response with a low dimensional projection of the predictors. More precisely,  SIR postulates the following {\it multiple index model}  for the data:
\begin{align}\label{model:multiple}
y=f(\bbeta_{1}^{\tau}\vx,\bbeta_{2}^{\tau}\vx,...,\bbeta_{d}^{\tau}\vx,\epsilon),
\end{align}
and estimates the subspace $\mathcal{S}=span\big\{~\bbeta_{1},...,\bbeta_{d}~\big\}$ via an eigen-analysis of the estimated conditional covariance matrix $var[\bbE(\vx | y)]$. Note that the individual $\bbeta_{i}$'s are not identifiable,  but the space  $\mathcal{S}$ can be estimated well. Based on the observation that $y\independent \vx \mid P_{\mathcal{S}}\vx$, 
\cite{cook1998regressiongraph} proposed a more general framework for dimension reduction without loss of information, 
often referred to as the {\it Sufficient Dimension Reduction} (SDR). 
Under this framework,  researchers look for the minimal subspace $\mathcal{S}'\subset \bbR^{p}$ such that $y\independent \vx \mid P_{\mathcal{S}'\vx}$ where $y$ is no longer necessarily a scalar response. Although numerous SDR algorithms have been developed in the past decades, SIR is still the most popular one among practitioners because of its simplicity and computational efficiency. Asymptotic theories developed for these SDR algorithms have all focused on scenarios where the data dimension $p$ is either fixed or growing at a much slower rate compared with the sample size $n$  \citep{dennis2000save,li2007directional,li2000high}.
The {\it `large p, small n'} characteristic of modern data raises  new challenges to these SDR algorithms. 

\cite{lin2015consistency} recently showed under mild conditions that the SIR estimate of the central space is consistent if and only if $\lim\frac{p}{n}=0$. This provides a theoretical justification for the necessity of the structural assumption such as sparsity for SIR when $p>n$. A commonly employed and also practically meaningful structural assumption made for high-dimensional linear regression problems is the sparsity assumption,  i.e., only a few predictors among the thousands or millions of candidate ones participate in the model.
We will show that this sparsity assumption can also  rescue {\it the curse of dimension} for dimension reduction algorithms such as SIR. 
Motivated by the Lasso and the regularized sparse PCA \citep{tibshirani1996regression,zou2005regularization},
\cite{li2006sparse} and  \cite{li2007sparse}  proposed some regularization approaches for SIR and SDR.  However, these  approaches often fail in  high dimensional numerical examples and are difficult to rectify because little is known about theoretical behaviors of these algorithms in high dimensional problems. The DT-SIR algorithm in \cite{lin2015consistency} and the Sparse-SIR algorithm in \cite{lin2016lasso}, however, have been shown to provide consistent estimations.
We agree with \cite{cook2012estimating} that a detailed understanding of ``the behaviour of these SDR estimators when $n$ is not large relative to $p$'' might be the key to efficient high-dimensional SDR algorithms. 
The main objective of the current paper is to understand the fundamental limits of the sparse SIR problem from a decision theoretic point of view. Such an investigation is not only interesting in its own right, but will also inform the development and evaluation of other SDR algorithms developed for high-dimensional problems.

\cite{neykov2015support} considered the (signed)-support recovery problem of the following class of single index models
\begin{align*}
y=f(\bbeta^{\tau}\vx,\epsilon) \quad \bbeta_{i} \in \{\pm 1/\sqrt{s}, 0\}, \quad  supp(\bbeta)=s 
\end{align*}
where $\vx \sim N(0,\vI_{p})$, $\epsilon\sim N(0,1)$.
Let $\Scale[.9]{\xi=\frac{n}{s\log(p)}}$, they proved that  {\bf 1)}  If  $\xi$ is sufficiently small, any algorithm fails to recover the (signed) support of $\bbeta$ with probability at least $1/2 $ and {\bf 2)} If $\xi$ is sufficiently large, the DT-SIR algorithm (see \cite{lin2015consistency} or Algorithm 1 below) can recover the (signed) support with probability  converging to 1 as $n \rightarrow \infty$.   
That is, the minimal sample size required to recover the support of $\bbeta$ is of order $s\log(p)$.
These results shed us some light on the possibility of obtaining the optimal rate of SIR-type algorithms in high dimension. 

SIR is widely considered as a `generalized eigenvector' problem \citep{chen1998can}. Inspired by recent advances in sparse PCA  \citep{amini2008high,johnstone2004sparse,cai2013sparse,birnbaum2013minimax,vu2012minimax},   where researchers aim at estimating the principal eigenvectors of the spiked models,  it is reasonable to expect a similar phase transition phenomenon \citep{johnstone2004sparse}, the signed support recovery \citep{amini2008high}, and the optimal rate \citep{cai2013sparse} for SIR when $\bSigma=\vI$. 
%The previous work, \cite{lin2015consistency} and \cite{neykov2015support} are mainly 
%In Sliced Inverse Regression, if $\bSigma=\vI_{p}$, the central space is spanned by the top $d$ eigenvectors of the matrix $var(\bbE[\vx|y])$. 
%Thus, if the loading vectors $\bbeta_{j}'$s are sparse, $var(\bbE[\vx|y])$ is `spiked' as well. 
%This similarity motivates us to investigate the phase transition phenomenon in \cite{lin2015consistency}. 
However, as it was pointed out in \cite{lin2015consistency}, the sample means in corresponding slices are neither independent nor identically distributed.
The usual concentration inequalities are not applicable.
This difficulty forced them to develop the corresponding deviation properties, i.e., the {\it `key lemma'} in \cite{lin2015consistency}.
On the other hand, the observation that the number $H$ of slices is allowed to be finite when $d$ is bounded (as we always require that $H>d$) suggests that a consistent estimate of the central space based on finite (e.g., $H$) sample means %$\overline{\vx}_{h,\cdot}$ 
is possible. 
This is again similar to the so-called {\it High dimensional low sample size }(HDLSS) scenario of PCA, which was first studied in \cite{jung2009pca} by estimating the principal eigenvectors based on finite samples. 
These connections suggest that theoretical issues in sparse SIR might be analogous to those in sparse PCA. However, our results in this article suggest that sparse linear regression is a more appropriate prototype for sparse SIR.
%To summarize, we should expect that the sparse PCA is the prototype of the sparse SIR problem  when $\bSigma=\vI_{p}$.

%\paragraph{Our contributions}  
The main contribution of this article is the determination of the minimax rate for  estimating the central space over two classes of models $\mathfrak{M}\left(p,d,\lambda,\kappa\right)$ and  $\mathfrak{M}_{s,q}(p,d,\lambda,\kappa)$, defined in \eqref{model:oracle} and \eqref{model:high:dim} respectively.
The risk of our interest is $\bbE[\|P_{\vV}-P_{\widehat{\vV}}\|_{F}^{2}]$, where $\vV$ is an orthogonal matrix formed by an orthonormal basis of $\mathcal{S}$, and $P_{\widehat{\vV}}$ is an estimate of $P_{\vV}$, the projection matrix associated with the orthogonal matrix $\vV$.
%The optimal rate of estimating the central space over $\mathfrak{M}\left(p,d,\lambda,\kappa\right)$, which we referred to as the {\it `Oracle risk'}, 
%is an extension  of the result in \cite{lin2015consistency} by including $\lambda_{\min}(var(\bbE[\vx|y]))$, the smallest non-zero eigenvalue of $var(\bbE[\vx|y])$, into consideration. 
%The upper bound of the optimal convergence rate over the sparse class $\mathfrak{M}_{s,q}(p,d,\lambda,\kappa)$ 
%(e.g. Theorem  \ref{thm:rsik:sparse:d:fixed:lambda:fixed} and \ref{thm:rsik:sparse:d=1}) follows from
%this {\it `Oracle risk' } and the technical tools developed in \cite{cai2013sparse}, \cite{gao2014minimax} and \cite{lin2015consistency}.
We construct an estimator (computationally unrealistic)
such that the risk of this estimator is of order $\frac{ds+s\log(ep/s)}{n\lambda}\wedge d$. 
Under mild conditions, we further demonstrate that the risk of any estimator is bounded below by $\frac{s\log(ep/s)}{n\lambda}\wedge 1$ if the dimension of the central space $d$ is bounded. 
Thus, the minimax rate of the risk  $\bbE[\|P_{\vV}-P_{\widehat{\vV}}\|_{F}^{2}]$ is $\frac{ds+s\log(ep/s)}{n\lambda}\wedge d$ if $d$ is bounded.
One of the key components of our analysis is the linear algebraic Lemma  \ref{lem:norm}, which might be of independent interest and be used in determining lower bounds of minimax estimation rates for other dimension reduction problems. 
%With this lemma, if we further assuming the technical conjecture \ref{conj:derivative} and that $d$ is bounded, we proved that the optimal rate of estimating the central space is $\frac{s\log(ep/s)}{n\lambda}$.
To the best of our knowledge, this is the first result about the minimax rate of estimating the central space in high dimension. %  including the information of eigenvalues.
In Subsection \ref{subsec:DT-SIR}, we show that the computationally efficient algorithm DT-SIR \citep{lin2015consistency}  achieves this optimal rate when $d=1$ and $s=O(p^{1-\delta})$ for some $\delta>0$. %which suggests us that a more appropriate prototype for high dimensional SIR ( and other SDR algorithms) might be sparse linear regression rather than generalized eigenvalue problems.
Furthermore, we investigate the effects of the slice number $H$ in the SIR procedure.

The rest of the paper is organized as follows. Section \ref{sec:main:results} presents the main results of the paper, including the rate of  the {\it oracle risk} in Section \ref{subsec:model:oracle} and the rate of the {\it sparse risk} in Section \ref{subsec:sparse:risk}. Since the lower bound can be obtained by  modifying some standard arguments, we defer its related proofs to the online supplementary file \citep{lin2016minimax} and give the
 proofs of upper bounds  in Sections \ref{sec:oracle:risk:proof} and  \ref{sec:sparse:risk:proof}. 
In Section \ref{sec:conclusion} we discuss potential extensions of our results. More auxiliary results and technical lemmas are included in the online supplementary file \citep{lin2016minimax}.

\section{Main Results} \label{sec:main:results} 
Since the establishment of the SDR framework about two decades ago,
estimating the central space has been investigated  under different assumptions   \citep{dennis2000save,cook1998regressiongraph,schott1994determining, ferre1998determining,li2007directional,hsing1992asymptotic,cook2012estimating}.
Various SDR algorithms have their own advantages and disadvantages for certain classes of link functions (models).
 For example, SIR only works when both the linearity and coverage conditions are satisfied \citep{li1991sliced}; 
Sliced Average Variance Estimation (SAVE) \citep{dennis2000save} works when the coverage condition is slightly violated but requires the constant variance condition. 
Thus, 
%it is formidable to ask for optimally estimating the central space of model \eqref{model:multiple} for arbitrary class of link functions.
to discuss the minimax rate of estimating the central space for model \eqref{model:multiple},  
it is necessary to first specify the class of models 
where one or several algorithms are practically used, 
and then check if these algorithms and their variants can estimate the central space optimally over this class of models.
SIR is one of the most well understood SDR algorithms, and is of special interests to know if it is rate optimal over a large class of models. 
This will not only improve our understanding of   high dimensional behaviors of SIR and its variants, but also  bring us  insights on  behaviors of other SDR algorithms.
%To avoid unnecessary confusion and complication, we assume the dimension $d$ of central space is bounded throughout this paper unless we explicitly stated (e.g., in Section \ref{subsec:upper_bound}).

\subsection{Notation}
In addition to those that have been used in Section 1, we adopt the following notations throughout the article.
For a matrix $\vV$, we denote its column space by $col(\vV)$ and its $i$-th row and $j$-th column by $\vV_{i,*}$ and $\vV_{*,j}$ respectively.
For vectors $\vx$ and $\bbeta$ $\in$ $\mathbb{R}^{p}$, we denote the $k$-th entry of $\vx$ as $\vx(k)$ and the inner product $\langle \vx, \bbeta \rangle$ as $\vx(\bbeta)$.
For two positive number a,b, we use $a\vee b$ and $a\wedge b$ to denote $\max\{a,b\}$ and  $\min\{a,b\}$, respectively.
For a matrix $A$, $\|A\|_{F}=tr(AA^{\tau})^{1/2}$. For a positive integer $p$, $[p]$ denotes the index set $\{1,2,...,p\}$. 
We use $C$, $C'$, $C_1$ and $C_2$ to denote generic absolute constants, though the actual value may vary from case to case.
For two sequences $a_{n}$ and $b_{n}$, we denote $a_{n}\succ b_{n}$ and $a_{n}\prec b_{n}$ if there exist positive constants $C$ and $C'$ such that $a_{n} \geq Cb_{n}$ and $a_{n} \leq C'b_{n}$, respectively. We denote $a_{n}\asymp b_{n}$ if both $a_{n}\succ b_{n}$ and $a_{n}\prec b_{n}$ hold.

\subsection{A brief review of SIR}\label{subsec:sir}
 Since we are interested in the space spanned by $\bbeta_{i}$'s in model \eqref{model:multiple}, without loss of generality, we can assume that $\vV=(\bbeta_{1},...,\bbeta_{d})$ is a $p\times d$ orthogonal matrix(i.e., $\vV^{\tau}\vV=\vI_{d}$) 
%The central subspace $\mathcal{S}$ is the minimal subspace $\mathcal{S}'$ such that $y\independent \vx\Big| P_{\mathcal{S}'}\vx$. It is easily seen that $\mathcal{S}$ is a subspace of the column space of $\vV$. 
and the models considered in this paper are
\begin{align}\label{model:modified:multiple}
y=f(\vV^{\tau}\vx,\epsilon), \quad \vV\in \mathbb{O}(p,d)
\end{align}
where $\vx \sim N(0,\vI_{p})$, $\epsilon\sim N(0,1)$, and $\mathbb{O}(p,d)$ is the set of all $p\times d$ orthogonal matrices.  Though $\vV$ is not identifiable, the column space $col(\vV)$ is estimable. 
The {\it Sliced Inverse Regression} (SIR) procedure proposed in \cite{li1991sliced}  estimate the central space $col(\vV)$ without knowing $f(\cdot)$, which can be briefly summarized as follows.
% SIR first estimates $\bLambda_{p}=var\left(\bbE[\vx|y]\right)$ via the sliced samples means. To be more precise,
Given $n$ $i.i.d.$ samples $(y_{i},\vx_{i})$, $i=1,\cdots,n$, SIR first divides them into $H$ equal-sized slices according to the order statistics $y_{(i)}$.\footnote{To ease notations and arguments, we assume that $n=cH$.}
We re-express the data as $y_{h,j}$ and $\vx_{h,j}$, where $(h,j)$ is the double subscript in which $h$ refers to the slice number and $j$ refers to the order number of a sample in the $h$-th slice, i.e.,
\begin{equation*}
y_{h,j}=y_{(c(h-1)+j)}, \mbox{\quad \quad } \vx_{h,j}=\vx_{(c(h-1)+j)}.
\end{equation*}
Here $\vx_{(k)}$ is the concomitant of $y_{(k)}$.
Let the sample mean in the $h$-th slice % $\frac{1}{c}\sum_{j=1}^c \vx_{h,j}$ 
be  $
\overline{\vx}_{h,\cdot}$, and the overall sample mean be $\vxbarbar$. 
SIR estimates $\bLambda\triangleq var(\bbE[\vx|y])$ by
\begin{equation}\label{eqn:lambda}
\widehat{\bLambda}_{H}=\frac{1}{H}\sum_{h=1}^{H}\bar{\vx}_{h,\cdot}\bar{\vx}_{h,\cdot}^{\tau}
\end{equation}
%and the SIR estimator $\widehat{\vV}_{H,c}$ of $\vV$ consists of the top $d$-eigenvectors of $\widehat{\bLambda}_{H,c}$.
and estimates the central space $col(\vV)$ by $col(\widehat{\vV}_{H})$  where $\widehat{\vV}_{H}$ is the matrix formed by the top $d$ eigenvectors 
 %$\widehat{\boldeta}^{\mathcal{I}}_1,\cdots, \widehat{\boldeta}^{\mathcal{I}}_d$ 
 of $\widehat{\bLambda}_{H}$.
Throughout this article, we assume that $d$, dimension of the central space, is known.

In order for the SIR to give a consistent estimate of the central space, following sufficient conditions have been suggested (e.g., \cite{li1991sliced}, \cite{hsing1992asymptotic} and \cite{zhu2006sliced}):
 \begin{itemize}
\item[${\bf A')}$]Linearity Condition and Coverage Condition:
\end{itemize}
\begin{align*}
span\Big\{~\bbE[\vx|y]~\Big\} = span\Big\{~\vV_{*,1},...,\vV_{*,d}~\Big\}
\end{align*}
where $\vV_{*,i}$ is the $i$-th columns of the orthogonal matrix $\vV$.
\begin{itemize}
\item[${\bf (B')}$] Smoothness and Tail conditions on the Central Curve $\bbE[\vx|y]$.
\end{itemize}
{\it Smoothness condition: } For $B>0$ and $n\geq 1$, let $\Pi_{n}(B)$ be the collection of all the $n$-point partitions $-B\leq y_{(1)} \leq \cdots \leq y_{(n)} \leq B$ of $[-B,B]$.  The central curve $\vm(y)$ satisfies the following conditions:
\begin{equation*}
\lim_{n\rightarrow \infty} \sup_{y\in \Pi_{n}(B)}n^{-1/4}\sum_{i=2}^{n}\|\vm(y_{i})-\vm(y_{i-1})\|_{2}=0, \forall B>0.
\end{equation*}
{\it Tail condition: } For some $B_{0}>0$, there exists a non-decreasing function $\widetilde{m}(y)$ on $(B_{0},\infty)$, such that
\begin{align}\label{cond:tail:hsing}
&\widetilde{m}^{4}(y)P(|Y|>y) \rightarrow 0 \mbox{ as } y\rightarrow \infty\\
\nonumber\|\vm(y)-\vm(y')\|_{2} &\leq |\widetilde{m}(y)-\widetilde{m}(y')| \mbox{ for } y, y' \in (-\infty,-B_{0})\cup(B_{0},\infty).
\end{align}

As in \cite{lin2015consistency}, where they demonstrated the phase transition phenomenon of SIR in high dimension,
we replace Condition {\bf (B$'$)} by 
\begin{itemize}
\item[${\bf (B'')}$] Modified Smoothness and Tail conditions,
\end{itemize}
which is all the same as  {\bf (B$'$)} except that eqn (\ref{cond:tail:hsing}) is replaced by
%{\it Smoothness condition: } For $B>0$ and $n\geq 1$, let $\Pi_{n}(B)$ be the collection of all the $n$-point partitions $-B\leq y_{(1)} \leq \cdots \leq y_{(n)} \leq B$ of $[-B,B]$.  The central curve $\vm(y)$ satisfies the following conditions:
%\begin{equation*}
%\lim_{n\rightarrow \infty} \sup_{y\in \Pi_{n}(B)}n^{-1/4}\sum_{i=2}^{n}\|\vm(y_{i})-\vm(y_{i-1})\|_{2}=0, \forall B>0.
%\end{equation*}
%{\it Tail condition: } For some $B_{0}>0$, there exists a non-decreasing function $\widetilde{m}(y)$ on $(B_{0},\infty)$, such that
\begin{align}\label{cond:stronger_tail}
& \bbE[\widetilde{m}(y)^{4}]<\infty\\
\nonumber\|\vm(y)-\vm(y')\|_{2} &\leq |\widetilde{m}(y)-\widetilde{m}(y')| \mbox{ for } y, y' \in (-\infty,-B_{0})\cup(B_{0},\infty).
\end{align}

%The conditions ${\bf A')}$ and ${\bf B')}$ in Subsection \ref{subsec:sir} are well accepted in literature.
%, however, it is unclear if we can determine the optimal rate over this class of models.  

%After replacing  the tail condition $\widetilde{m}^{4}(y)P(|Y|>y) \rightarrow 0 \mbox{ as } y\rightarrow \infty$ in {\bf B$'$)}  by $\bbE[\widetilde{m}(y)^{4}]<\infty$ in {\bf B$''$}. 
 
It is easy to see that Condition {\bf (B$''$)} is slightly stronger than Condition {\bf (B$'$)}. A main advantage of Condition {\bf (B$''$)} is the following proposition proved   in \cite{neykov2015support}.
\begin{prop}\label{prop:sliced}
%If we replace the tail condition, $\widetilde{m}^{4}(y)P(|Y|>y) \rightarrow 0 \mbox{ as } y\rightarrow \infty$,  in \eqref{cond:tail:hsing} by a slightly stronger condition
%\begin{align*}
%\bbE[\widetilde{m}(Y)^{4}]<\infty,
%\end{align*}
If Condition {\bf B$''$} holds, the central curve $\bbE[\vx|y]$ satisfies the sliced stable condition  (defined below) with $\vartheta=\frac{1}{2}$.
\end{prop} 
 
 \begin{definition}\label{sliced-stable}
 Let $Y$ be a random variable. For  $0<\bgamma_{1} < 1<\bgamma_{2} $, 
 let $\mathcal{A}_H(\bgamma_{1},\bgamma_{2})$ denote all partitions $\{ -\infty = a_0 \leq a_2 \leq \ldots \leq a_{H} = +\infty\}$ of $\mathbb{R}$, such that 
 \[
 \frac{\bgamma_{1}}{H} \leq \mathbb{P}(a_h \leq Y \leq a_{h + 1}) \leq \frac{\bgamma_{2}}{H}.
 \]
A curve $\vm(y)$ is $\vartheta$-sliced stable with respect to Y, if there exist positive constants $\bgamma_{1},\bgamma_{2},\bgamma_{3}$ such that for any partition $\in \mathcal{A}_H(\bgamma_{1},\bgamma_{2})$ and any $\bbeta \in \bbR^{p}$ , we have
\begin{align}\label{eqn:temp:sliced}
\frac{1}{H}\sum_{h=1}^{H}var\left(\bbeta^{\tau}\vm(Y)\big| a_{h-1}\leq Y < a_{h}\right)\leq \frac{\bgamma_{3}}{H^{\vartheta}}  var\left(\bbeta^{\tau}\vm(Y)\right).
\end{align}
A curve is sliced stable if it is $\vartheta$-sliced stable for some positive constant $\vartheta$.
 \end{definition}

%\begin{definition}\label{def:sliced:stable} 

%\end{definition}

Intuitively,  $H$ $\rightarrow \infty$ implies that the LHS of \eqref{eqn:temp:sliced} 
converges to zero.
%$\rightarrow var(\bbE[\bbeta^{\tau}\vm(y)|y])=0$. 
%In order to utilize this approximation,  a quantitative description of this convergence rate is preferred. 
Definition~\ref{sliced-stable} states that its convergence rate  is a power of $H$, 
although any function of $H$ that converges to 0 can be placed before $var(\bbeta^{\tau}\vx)$  on the RHS of \eqref{eqn:temp:sliced}. 
Thus, the sliced stable condition is almost the necessary condition to ensure that the SIR works. A main advantage of the sliced stable condition is that we can easily quantify the deviation properties of the eigenvalues, eigenvectors, and each entries of $\widehat{\bLambda}_{H}$.  
This is one of the main technical contributions of \cite{lin2015consistency}.
%For example, it is easy to see that  for each slice interval $S_{h}=(y_{h,c},y_{h+1,c}]$, we have $\bbP(S_{h}) \approx \frac{1}{H}$ with high probability. 
%Conditioning on this event, \cite{lin2015consistency} proved their {\it `key lemma' }  in (cf. Lemma \ref{lem:main:deviation} in this paper). 
We henceforth assume that the central curve satisfies the sliced stable condition.
As shown by  Proposition  \ref{prop:sliced}, Condition {\bf (B$''$)} ensures the sliced-stable condition.

\subsection{The class of functions $\mathcal{F}_{d}(\lambda ,\kappa)$}
Let $\vz=\vV^{\tau}\vx$, then $\vz \sim N(0,\vI_{d})$. Let $\bLambda_{\vz}=var(\bbE[\vz|y])$. 
Since 
$
\bbE[\vx|y]=P_{\vV}\bbE[\vx|y]=\vV\bbE[\vV^{\tau}\vx|y]=\vV\bbE[\vz|y]$,  the sliced stability for $\bbE[\vz|y]$ implies the sliced stability for $\bbE[\vx|y]$
and vice verse.  
Since we have assumed that $\vx\sim N(0,\vI_{p})$, the linearity condition holds automatically. The coverage condition, which requires $rank(var(\bbE[\vx|y]))=d$, can be refined as
\begin{align}
\lambda\leq\lambda_{d}(var(\bbE[\vx|y])) \leq \lambda_{1}(var(\bbE[\vx|y]))\leq \kappa\lambda \leq 1
\end{align}
for some positive constant $\kappa>1$.
Since $\bLambda \triangleq  var(\bbE[\vx|y])=\vV\bLambda_{\vz}\vV^{\tau}
$,
we know $\lambda_{j}(\bLambda)=\lambda_{j}(\bLambda_{\vz}), j=1,...,d$. In particular, we have
$
\lambda\leq\lambda_{d}(var(\bbE[\vz|y])) \leq \lambda_{1}(var(\bbE[\vz|y]))\leq \kappa\lambda \leq 1,$  where $\kappa$ is assumed to be a fixed constant. This coverage condition is commonly adopted in the literature (e.g., \cite{cai2013sparse} and \cite{gao2014minimax}) when researchers discuss the dimension reduction problems.
The class of functions $f$ satisfying the sliced stable condition and coverage condition is of our main interests in this paper. More precisely, we introduce $\mathcal{F}_{d}(\lambda,\kappa)$ as below.

\begin{definition}
Let $\vz \sim N(0,\vI_{d})$ and $\epsilon \sim N(0,1)$. 
A  function $f(\vz,\epsilon)$ belongs to the class $\mathcal{F}_{d}(\lambda ,\kappa)$, if  the following conditions are satisfied.
\end{definition}
\begin{itemize}
\item[${\bf (A)}$] Coverage condition:
$0<\lambda \leq \lambda_{d}(\bLambda_{z}) \leq ...\leq \lambda_{1}(\bLambda_{z}) \leq  \kappa\lambda \leq 1, $
where  $\bLambda_{z}\triangleq var(\bbE[\vz|f(\vz,\epsilon)])$.

\item[${\bf (B)}$] Sliced stable condition:
$\vm_{z}(y)=\bbE[\vz| f(\vz, \epsilon)]$  is sliced stable with respect to $y$, 
where  $y=f(z,\epsilon)$.
\end{itemize}

It is easy to see that almost all functions $f$  that make SIR work belong to $\mathcal{F}_{d}(\lambda,\kappa)$ for some $\kappa$ and $\lambda$.

\subsection{Upper bound of the risk}\label{subsec:upper_bound}
Suppose  we have $n$ samples generated from a multiple index model $\mathcal{M}$ with link function $f$ and orthogonal matrix $\vV$, that is, $y=f(\vV^{\tau}\vx,\epsilon)$.  
We are interested in the risk $\bbE_{\mathcal{M}}\|P_{\widehat{\vV}}-P_{\vV}\|^{2}_{F}$ where $P_{\widehat{\vV}}$ is an estimate of $P_{\vV}$ based on these samples. 
In this subsection, we provide an upper bound of this risk.
All detailed proofs are deferred to Sections \ref{sec:oracle:risk:proof},  \ref{sec:sparse:risk:proof}, and online supplementary file \citep{lin2016minimax}.

\subsubsection{Oracle Risk} \label{subsec:model:oracle} Here we are interested in estimating the central space over the following class of models parametrized by $(\vV,f)$:
\begin{align} \label{model:oracle}
\mathfrak{M}\left(p,d,\lambda,\kappa\right)\triangleq \Big\{~(\vV, f) ~\Big|~ \vV \in \mathbb{O}(p,d), f \in \mathcal{F}_{d}(\lambda ,\kappa) ~\Big\}.
\end{align}
A main result of this article is:
\begin{thm}[An Upper Bound of Oracle Risk] 
Assume that $\frac{dp}{n\lambda}$ is sufficiently small and $d^{2}\leq p$. We have
\label{thm:risk:oracle:upper:d}
\begin{align}\label{eqn:oracle:risk}
\inf_{\widehat{\vV}} \sup_{\mathcal{M} \in \mathfrak{M}\left(p,d,\lambda,\kappa\right)} \bbE_{\mathcal{M}}\|P_{\widehat{\vV}}-P_{\vV}\|^{2}_{F}\prec d\wedge\frac{d(p-d)}{n\lambda}.
\end{align}
\end{thm}  
In order to establish the upper bound, we consider the estimate
$\widehat{\vV}_{H}$, which is a $p\times d$ orthogonal matrix forming by the top-$d$ eigenvectors of $\widehat{\bLambda}_{H}$,
and show that $P_{\widehat{\vV}_{H}}$ achieves the rate in Theorem \ref{thm:risk:oracle:upper:d}.

A result in \cite{lin2015consistency}, which states that 
\begin{align}\label{eqn:OLD:results}
\|\widehat{\bLambda}_{H}-var(\bbE[\vx|y])\|_{2}=O_{P}\left(
\frac{1}{H^{\vartheta}}+\frac{H^{2}p}{n}+\sqrt{\frac{H^{2}p}{n}}\right),
\end{align}
appears to contradict our Theorem \ref{thm:risk:oracle:upper:d} here: 
(i) it does not depend on $d$, the dimension of central subspace;
(ii) it does not depend on $\lambda$, the smallest non-zero eigenvalue of $var(\bbE[\vx|y])$;
(iii) it depends on $H$ (the number of slices) and seems worse than our upper bound here. 
The first two differences appear simply because \cite{lin2015consistency} have assumed that $d$ is bounded and the non-zero eigenvalues of $var(\bbE[\vx|y])$ are bounded below by some positive constant (i.e., the information about eigenvalues and $d$ is absorbed by some constants). The third difference appears
because  we here are interested in the convergence rate of the SIR estimate of the  space $\mathcal{S}$  rather than the convergence rate of the SIR estimate of the matrix $var(\bbE[\vx|y])$.
As they have pointed out, the convergence rate of $\widehat{\bLambda}_{H}$ might be different (slower) than the convergence rate of $P_{\widehat{\vV}_{H}}$. More precisely, we have
\begin{align}\label{eqn:temp:temp:rmk}
\widehat{\bLambda}_{H}-\bLambda = \left( \widehat{\bLambda}_{H}-P_{\vV}\widehat{\bLambda}_{H}P_{\vV}\right)
+\left(P_{\vV}\widehat{\bLambda}_{H}P_{\vV}-\bLambda\right).
\end{align}
From the proof of Theorem 1 of \cite{lin2015consistency}, we can easily check that
 the first term is of rate $\frac{pH^{2}}{n}+\sqrt{\frac{pH^{2}}{n}}$ and the second term is of  rate $\frac{1}{H^{\vartheta}}$.  
Since $P_{\vV}\widehat{\bLambda}_{H}P_{\vV}$ and $\bLambda$ share the same column space and we are interested in estimating $P_{\vV}$, the convergence rate of the second term in \eqref{eqn:temp:temp:rmk} does not matter provided that $H$ is a large enough integer.
Thus, Theorem \ref{thm:risk:oracle:upper:d} does not contradict the convergence result in \cite{lin2015consistency}.

\begin{rmk} On the role of $H$. \normalfont Researchers have claimed that the performance of SIR procedure is not sensitive to the choice of $H$, i.e., $H$ can be as large as $\frac{n}{2}$ \citep{hsing1992asymptotic} and can also be a large enough fixed integer when $d=1$ \citep{duan1991slicing}. A direct corollary of  Theorem \ref{thm:risk:oracle:upper:d} is that if $d$ is fixed, $H$ can be a large enough constant such that $col(\widehat{\vV}_{H})$ is an optimal estimate of $col(\vV)$. 
In the SIR literature, researchers care about the eigenvectors of $\bLambda$ and ignore the eigenvalue information. 
In this article, we show that  the larger the $H$, the more accurate the estimate of the eigenvalues of $\bLambda$, and  
illustrate this phenomenon via numerical simulations in Section \ref{subsec:numerical:H}.
Taking the eigenvalue information into consideration will bring us more a detailed understanding of SIR.

\end{rmk}

\subsubsection{Upper bound of the risk of sparse SIR}\label{subsec:sparse:risk}
%Theorem \ref{thm:risk:oracle:upper:d} 
\cite{lin2015consistency} shows that when dimension $p$ is larger than or comparable with the sample size $n$, the SIR estimate of the central space  is inconsistent. Thus, structural assumptions such as sparsity are necessary for high dimensional SIR problem.We here impose the weak $l_{q}$ sparsity on the loading vectors $\vV_{*,1},...,\vV_{*,d}$.
For a $p\times d$ orthogonal matrix $\vV$ (i.e., $\vV^{\tau}\vV=\vI_{d}$), we order the row norms in decreasing order as $\|\vV_{(1),*}\|_{2}\geq ... \geq \|\vV_{(p),*}\|_{2}$ and define the weak $l_{q}$ radius of $\vV$ to be
\begin{align}\label{def:sparse orthogonal}
\|\vV\|_{q,w}\triangleq \max_{j\in [p]} j\|\vV_{(j),*}\|^{q}.
\end{align}
Let
$\mathbb{O}_{s,q}(p,d)=\big\{ \vV ~\big| ~ \vV \in \mathbb{O}(p, d) \mbox{ such that }  \|\vV\|_{q,w} \leq s ~\big\}
$ be the set of weak $l_{q}$ sparse orthogonal matrices.
Weak $l_{q}$-ball is a commonly used condition for sparsity. See, for
example, \cite{abramovich2006special} for wavelet estimation and \cite{cai2012optimal} for sparse co-variance matrix estimation.
Furthermore, we need the notion of {\it effective support }, which was introduced  by \cite{cai2013sparse}. 
The size of {\it effective support } is defined to be
$
k_{q,s}\triangleq \lceil x_{q}(s,d) \rceil,
$
where 
\begin{align}\label{definition of k}
x_{q}(s,d)\triangleq \max\Big\{ 0\leq x \leq p~ |~ x \leq s \left(\frac{n\lambda}{d+\log\left(\frac{ep}{x} \right)} \right)^{q/2} \Big\}
\end{align}
and $\lceil a \rceil$ denotes the smallest integer no less than $a \in \mathbb{R}$.
For more detailed discussions of the sparse orthogonal matrices, we refer to \cite{cai2013sparse}.

In this subsection, we are interested in estimating the central space over the following class of high dimensional models parametrized by $(\vV,f)$:
\begin{align}\label{model:high:dim}
\mathfrak{M}_{s,q}\left( p,d,\lambda,\kappa\right) \triangleq \Big\{(\vV, f) ~\Big| ~\vV \in \mathbb{O}_{s,q}(p,d), f \in \mathcal{F}_{d}(\lambda ,\kappa) \Big\}. 
\end{align}
Let 
$
\epsilon_{n}^{2} \triangleq \frac{1}{n\lambda}\left(dk_{q,s}+k_{q,s}\log\frac{ep}{k_{q,s}} \right)$. 
We have the following result:  
\begin{thm}[The Upper Bound of Optimal Rates]
\label{thm:risk:sparse:upper:d}
Assume that $\kappa$ is fixed, $d^{2}\leq k_{q,s}$, $\epsilon^{2}_{n}$ is sufficiently small and $n\lambda \leq e^{p}$. We have
\begin{align}\label{eqn:rate:sparse:rsik:s=q}
\inf_{\widehat{\vV}} \sup_{\mathcal{M}\in \mathfrak{M}_{s,q}\left( p,d,\lambda,\kappa\right) } \bbE_{\mathcal{M}} \|P_{\widehat{\vV}}-P_{\vV}\|^{2}_{F}\prec d \wedge \frac{dk_{q,s}+k_{q,s}\log\frac{ep}{k_{q,s}} }{n\lambda}.
\end{align}
\end{thm}
In order to establish the upper bound in Theorem \ref{thm:risk:sparse:upper:d}, we need to construct an estimator that attains it.
Let $\mathcal{B}(k_{q,s})$ be the set of all subsets of $[p]$ with size $k_{q,s}$. To ease the notation,  we often drop the subscript $(q,s)$ of $k_{q,s}$ below and assume that there are $n=2Hc$ samples. 
Let us divide the samples randomly into two equal size sets. 
Let $\widehat{\bLambda}^{(1)}_{H}$ and $\widehat{\bLambda}^{(2)}_{H}$ be the SIR estimates of $\bLambda=var(\bbE[\vx|y])$ based on the first and second sets of samples, respectively.  
Inspired by the idea in  \cite{cai2013sparse}, we introduce the following aggregation estimator $\widehat{\vV}_{E}$ of $\vV$.
\paragraph{\it Aggregation Estimator $\widehat{\vV}_{E}$}
For each $B \in \mathcal{B}_{k}$, we let 
 \begin{equation}\label{estimator:sample}
\begin{aligned}
&\widehat{\vV}_{B}\triangleq \arg\max_{\vV} \langle\widehat{\bLambda}_{H}^{(1)},\vV\vV^{\tau} \rangle=\arg\max_{\vV}Tr(\vV^{\tau}
\widehat{\bLambda}_{H}^{(1)}\vV)\\
&\mbox{ s.t. }  \vV^{\tau}\vV=\vI_{d},   \|\vV\|_{q,w}=k \mbox{ and } supp(\widehat{\vV}_{B})\subset B
\end{aligned}
\end{equation}
and
\begin{align*}
B^{*}\triangleq \arg \max_{B\in \mathcal{B}(k)} \langle\widehat{\bLambda}_{H}^{(2)},\widehat{\vV}_{B}\widehat{\vV}_{B}^{\tau} \rangle=\arg \max_{B\in \mathcal{B}(k)}Tr(\widehat{\vV}_{B}^{\tau}
\widehat{\bLambda}_{H}^{(2)}\widehat{\vV}_{B}).
\end{align*}
Our aggregation estimator $\widehat{\vV}_{E}$ is defined to be $\vV_{B*}$. 

\vspace*{3mm}

$B^{*}$ is a stochastic set and, for any fixed $B$, $\widehat{\vV}_{B}$ is independent of the second set of samples.
From the definition of $\widehat{\vV}_{E}$, it is easy to see
\begin{align}\label{eqn:two:set:inequality}
\langle\bLambda^{(2)}_{H},\widehat{\vV}_{E}\widehat{\vV}^{\tau}_{E}-\widehat{\vV}_{B}\widehat{\vV}^{\tau}_{B} \rangle \geq 0
\end{align}
for any $\widehat{\vV}_{B}$ where $B \in \mathcal{B}$. In Section \ref{sec:sparse:risk:proof}, we will show that the aggregation estimator $\widehat{\vV}_{E}$ achieves the converges rate on the right hand side of \eqref{eqn:rate:sparse:rsik:s=q}.

\subsection{Lower Bound and Minimax Risk} 
To avoid unnecessary details,  we assume that  dimension $d$ of the central space is bounded in this subsection.
%Comparing with other parametric models ( e.g., linear models, spiked models for sparse PCA), 
The semi-parametric characteristic of multiple index models brings us additional difficulties in determining the lower bound of the minimax rate. 
Because of our ignorance on the function class $\mathcal{F}_{d}(\lambda,\kappa)$, we can only establish the lower bound in two restrictive  cases: 
(i) $\lambda$, the smallest non-zero eigenvalue of $var(\bbE[\vx|y])$, is bounded below by some positive constant; 
and (ii) single index models where $d=1$.
%Although the convergence rate depending optimally on all the parameters $n$, $d$, $\lambda$, $s$ and  $p$ is more desired, 
%we feel that our partial results are worth to report.  
To the best of our knowledge, even the optimal rate of estimating the central space depending only on $n$, $s$ and $p$ in high dimensions has never been discussed in the literature.
Furthermore, we have observed from extensive numerical studies that the $4$-th direction is difficult to detect for $p=10$ even with the sample size greater than $10^{6}$. 
This observation conforms to the existing numerical studies reported in the literature, i.e., most researchers only reported numerical studies for models with $d\leq 2$ except that \cite{ferre1998determining} performed a numerical study for a model with $d=4$ and reported that the $4$-th direction was hard to discover. Thus, the optimal rate with $d$ bounded  might be a more reasonable target to pursue.

\subsubsection{ $\lambda$ is bounded below by some positive constant} Assume that $\lambda$, 
the smallest non-zero eigenvalues of $var(\bbE[\vx|y])$, 
is bounded below by a positive constant. 
We have the following optimal convergence rate of the {\it Oracle Risk}.
\begin{thm}[ Oracle Risk] \label{thm:oracle:risk:d:fxied:lambda:fixed}Assume that $d, \lambda$ are bounded. We have 
\begin{align}\label{eqn:oracle:risk:d:fixed}
\inf_{\widehat{\vV}} \sup_{\mathcal{M} \in \mathfrak{M}\left(p,d,\lambda,\kappa\right)} \bbE_{\mathcal{M}}\|P_{\widehat{\vV}}-P_{\vV}\|^{2}_{F}\asymp d\wedge\frac{dp}{n}.
\end{align}
\end{thm}
\begin{rmk}\normalfont
Although we have assumed that the dimension of the central space $d$ is bounded, we include it in the convergence rate to emphasize that the result holds for  multiple index models.
\end{rmk}
Because of Theorem \ref{thm:risk:oracle:upper:d},  we only need to  establish the lower bound.
We defer the detailed proof to the online supplementary file \citep{lin2016minimax} and briefly sketch its key steps here. 
One of the key steps in obtaining the lower bound is constructing a finite family of distributions that are distant from each other in the parameter space and close to each other in terms of the KL-divergence. Recall that, for any sufficiently small $\epsilon>0$ and any positive constant $\alpha<1$, \cite{cai2013sparse} have constructed a subset $\Theta \subset$ $\mathbb{G}(p,d)$, the Grassmannian manifold consisting of all the $d$ dimensional subspaces in $\bbR^{p}$, such that 
\begin{align*}
\left| \Theta\right|&\geq \left(\frac{c_{0}}{\alpha c_{1}}\right)^{d(p-d)} \mbox{ and }\\
\alpha^{2}\epsilon^{2}\leq \|\theta_{i}&-\theta_{j}\|_{F}^{2}\leq \epsilon^{2} \mbox{ for any } \theta_{i}, \theta_{j} \in \Theta
\end{align*}
 for some absolute constants $c_{0}$ and $c_{1}$. For any $\theta_{j} \in \Theta$,  if we can choose a $p\times d$ orthogonal matrix $\vB_{j}$ such that the column space of $\vB_{j}$ corresponds to $\theta_{j}\in \mathbb{G}(p,d)$,  we may  consider the following finite class of models
\begin{align*}
y=f(\vB_{j}^{\tau}\vx)+\epsilon, \vx\sim N(0,\vI_{p}) \mbox{ and } \epsilon \sim N(0,1)
\end{align*}
where $f$ is a $d$-variates function with bounded first derivative such that these models belong to  $ \mathfrak{M}\left(p,d,\lambda,\kappa\right)$. Let $p_{f,\mathbf{B}}$ denote the joint density of $(y,\vx)$. Simple calculation shows that
\begin{align}\label{eqn:kl:bound}
KL(p_{f,\vB_{1}},p_{f,\vB_{2}})\leq C\|\nabla f\|^{2}\|\vB_{1}-\vB_{2}\|^{2}_{F}\leq C \|\vB_{1}-\vB_{2}\|^{2}_{F}.
\end{align}
If we have
\begin{align}\label{eqn:distance:bound}
\|\vB_{1}-\vB_{2}\|^{2}_{F} \leq \|P_{\vB_{1}}-P_{\vB_{2}}\|^{2}_{F},
\end{align}
we may apply the standard Fano type argument ( e.g., \cite{cai2013sparse} ) to obtain the essential rate $\frac{dp}{n}$ of the lower bound.  

However, \eqref{eqn:distance:bound} is not always true (e.g., it fails if $\vB_{1}$ and $\vB_{2}$ are two different orthogonal matrices sharing the same column space).
We need to carefully specify  $\vB_{j}$ for each $\theta_{j} \in \Theta \subset$ $\mathbb{G}(p,d)$ such that they satisfy the inequality \eqref{eqn:distance:bound}. 
It seems to be a simple linear algebraic problem, however, its proof  requires (slightly) non-trivial work in differential geometry (cf. Lemma \ref{lem:norm}).
\footnote{Q. Lin appreciates the helpful discussions with Dr. Long Jin}  
Thus we know that the rate in Theorem \ref{thm:risk:oracle:upper:d} is optimal if $d$ and $\lambda$ are bounded.
Once the `{\it Oracle risk'} has been established, the standard argument in \cite{cai2013sparse} leads us the following:
\begin{thm}[Optimal Rates ]\label{thm:rsik:sparse:d:fixed:lambda:fixed}
Assume that $d$, $\lambda$ are bounded, and $n\lambda\leq e^{p}$. We have
\begin{align}\label{eqn:rate:sparse:rsik:s=q:d:fixed:lambda}
\inf_{\widehat{\vV}} \sup_{\mathcal{M}\in \mathfrak{M}_{s,q}\left( p,d,\lambda,\kappa\right) } \bbE_{\mathcal{M}} \|\widehat{\vV}\widehat{\vV}^{\tau}-\vV\vV^{\tau}\|^{2}_{F}\asymp d\wedge \frac{dk_{q,s}+k_{q,s}\log\frac{ep}{k_{q,s}} }{n}
\end{align}
\proof See the online supplementary file \citep{lin2016minimax}.
\end{thm}

%\begin{rmk}\normalfont
%In Theorem \ref{thm:rsik:sparse:d:fixed:lambda:fixed}, since $d$ is bounded, the above rate is equivalent to $(k_{q,s}/n)\log(ep/k_{q,s})$. We state it in current form in order to emphasize that $d$ is not necessarily to be 1.
%\end{rmk}

\subsubsection{Single Index Models}\label{subsub:d=1}
If we restrict our consideration to single index models (i.e., $d=1$), we have a convergence rate optimally depending on $n$, $\lambda$, $s$, and $p$.

\begin{thm}[Oracle Risk for Single Index Models] Assuming that $d=1$ and $n\lambda\leq e^{p}$, we have \label{thm:risk:oracle:d=1}
\begin{align}\label{eqn:oracle:risk:d:fixed}
\inf_{\widehat{\vV}} \sup_{\mathcal{M} \in \mathfrak{M}\left(p,d,\lambda,\kappa\right)} \bbE_{\mathcal{M}}\|\widehat{\vV}\widehat{\vV}^{\tau}-\vV\vV^{\tau}\|^{2}_{F}\asymp 1\wedge\frac{p}{n\lambda}.
\end{align}
\end{thm}
Since we have proved Theorem \ref{thm:risk:sparse:upper:d}, all we need to do is to establish a suitable lower bound. Let us consider the following linear model:
 \[
 y=f_{\lambda}(\bbeta^{\tau}\vx)=\sqrt{\lambda}\bbeta^{\tau}\vx+\epsilon,
 \]
  where $\bbeta$ is a unit vector, $\vx\sim N(0,\vI)$ and $\epsilon\sim N(0,1)$.    
Simple calculation shows that
\[
var(\bbE[\vx|y])=\frac{\lambda}{1+\lambda}\mbox{ and } |\nabla f_{\lambda}|\leq C\sqrt{\lambda}.
\] 
Thus, inequality \eqref{eqn:kl:bound} becomes
\begin{align}\label{eqn:kl:bound:lambda}
KL(p_{f,\bbeta_{1}},p_{f,\bbeta_{2}})\leq C\|\nabla f\|^{2}\|\bbeta_{1}-\bbeta_{2}\|^{2}_{F}\leq C\lambda \|\bbeta_{1}-\bbeta_{2}\|^{2}_{F}
\end{align}
and the desired lower bound follows from the same argument as that of Theorem  \ref{thm:oracle:risk:d:fxied:lambda:fixed}.
Once the oracle risk has been established, the standard argument in \cite{cai2013sparse} leads us to the following result:
 
\begin{thm}[Optimal Rates : $d=1$  ]\label{thm:rsik:sparse:d=1}
Assume that $d=1$ and  $n\lambda\leq e^{p}$. We have
\begin{align}\label{eqn:rate:sparse:rsik:s=q:d:fixed}
\inf_{\widehat{\vV}} \sup_{\mathcal{M}\in \mathfrak{M}_{s,q}\left( p,d,\lambda,\kappa\right) } \bbE_{\mathcal{M}} \|\widehat{\vV}\widehat{\vV}^{\tau}-\vV\vV^{\tau}\|^{2}_{F}\asymp 1\wedge \frac{k_{q,s}\log\frac{ep}{k_{q,s}} }{n\lambda}.
\end{align}
\proof It is similar to the proof of Theorem \ref{thm:rsik:sparse:d:fixed:lambda:fixed} and thus omitted. 
\end{thm}

\subsubsection{Multiple Index Models with $d$ bounded}
The arguments in the subsection \ref{subsub:d=1} motivate us to propose the following (conjectural) property for the function class $\mathcal{F}_{d}(\lambda,\kappa)$. 
\begin{conj} \label{conj:derivative} If $d$ is bounded, there is a constant $C$ such that for any $0<\lambda\leq 1$, there exists a $d$-variate function $f_{\lambda}$ such that $f_{\lambda}(x_{1},...,x_{d})+x_{d+1} \in \mathcal{F}_{d}(\lambda,\kappa)$ and
\begin{align}\label{conjecture:require}
\|\nabla f_{\lambda}(x_{1},...,x_{d})\|\leq C\sqrt{\lambda}. 
\end{align}
\end{conj}
\begin{rmk}\normalfont
Inequality \eqref{conjecture:require} can be relaxed to that $\|\nabla f_{\lambda}(x)\|\leq C\sqrt{\lambda}$ holds with high probability when $x \sim N(0,\vI_{d})$. 
\end{rmk}
The construction in subsection \ref{subsub:d=1} shows that this conjecture holds for $d=1$. For any $d>1$, suppose that there exists a function $f$  such that $f(x_{1},...,x_{d})+x_{d+1}\in \mathcal{F}_{d}(\mu,\kappa)$. We expect that, for $y=\sqrt{\lambda}f(\vx)+\epsilon$,  there exist constants $C_{1}$ and $C_{2}$ such that
\[
 C_{1}\lambda\leq \lambda_{d}(var(\bbE_{\lambda}[\vx|y]))\leq \lambda_{1}(var(\bbE_{\lambda}[\vx|y])) \leq  C_{2}\kappa\lambda.
\]
Note that the density function $p(y)$ of $y$ is the convolution of the density functions of $\epsilon$ and $\sqrt{\lambda}f(\vx)$. 
Heuristically, if $f(\vx)$ is (nearly) normal, by the continuity of the convolution operator, we expect that  $\lambda_{d}(var(\bbE[\vx|y]))$ $\asymp\lambda$.
Since we cannot prove it rigorously,  
we present some supporting numerical evidences here in Subsection \ref{subsubsection:numerical:conjecture}. 
Assuming this conjecture, we have the following theorems, of which the proofs are similar to those of Theorem \ref{thm:oracle:risk:d:fxied:lambda:fixed} and Theorem \ref{thm:rsik:sparse:d:fixed:lambda:fixed}.
\begin{thm}[ Oracle Risk : $d$ is bounded] Assuming that $d$ is bounded and  Conjecture \ref{conj:derivative} holds, we have \label{thm:oracle:risk}
\begin{align}\label{eqn:oracle:risk:d:fixed}
\inf_{\widehat{\vV}} \sup_{\mathcal{M} \in \mathfrak{M}\left(p,d,\lambda,\kappa\right)} \bbE_{\mathcal{M}}\|\widehat{\vV}\widehat{\vV}^{\tau}-\vV\vV^{\tau}\|^{2}_{F}\asymp d\wedge\frac{dp}{n\lambda}.
\end{align}
\proof It is similar to the proof of Theorem \ref{thm:oracle:risk:d:fxied:lambda:fixed}, and thus omitted. 
\end{thm}
\begin{thm}[Optimal Rates : $d$ is bounded ]\label{thm:sparse:rsik}
Assuming that $d$ is fixed, $n\lambda\leq e^{p}$  and Conjecture \ref{conj:derivative}  holds, we have
\begin{align}\label{eqn:rate:sparse:rsik:s=q:d:fixed}
\inf_{\widehat{\vV}} \sup_{\mathcal{M}\in \mathfrak{M}_{s,q}\left( p,d,\lambda,\kappa\right) } \bbE_{\mathcal{M}} \|\widehat{\vV}\widehat{\vV}^{\tau}-\vV\vV^{\tau}\|^{2}_{F}\asymp d\wedge \frac{dk_{q,s}+k_{q,s}\log\frac{ep}{k_{q,s}} }{n\lambda}.
\end{align}
\proof It is similar to the proof of Theorem \ref{thm:rsik:sparse:d:fixed:lambda:fixed}, and thus omitted. 
\end{thm}

\subsection{Optimality of DT-SIR} \label{subsec:DT-SIR}
In the previous section, we have proved that 
the aggregation estimator $\widehat{\vV}_{E}$ is rate optimal. 
In practice, however, it is computationally too expensive.
%We are interested in  if any computationally efficient algorithm is rate optimal.
The DT-SIR algorithm proposed in \cite{lin2015consistency} is computationally efficient in general, and can be further simplified when $\bSigma_{\vx}=\vI$. 

\begin{algorithm}
\caption{DT-SIR}
\begin{algorithmic}[1]\label{alg:DTSIR}
\vspace*{1mm}
\\ Let $S=\{~ i~ |~ \widehat{\bLambda}_{H}(i,i)>t ~\}$ for a properly choosen $t$.
\vspace*{1mm}
\\ Let $\widehat{\bbeta}$ be the principal eigenvector of $\widehat{\bLambda}_{H}(S,S)$.
\vspace*{1mm}
\\ We embed $\widehat{\bbeta}$ into $\mathbb{R}^{p}$ by filling the entries outside $S_{t}$ with 0 and denote it by $\widehat{\bbeta}_{DT}$.
\end{algorithmic}
\end{algorithm}

In this section, we focus on the single index model. 
with the exact sparsity on the loading vector $\bbeta$, i.e., $|supp(\bbeta)|=s$.
\begin{thm}\label{thm:DTSIR}
Suppose  $s=O(p^{1-\delta})$ for some $\delta>0$, $\frac{s\log(p)}{n\lambda}$ is sufficiently small and $n=O(p^{C})$ for some constant $C$.
Let $\widehat{\bbeta}_{DT}$ be the DT-SIR estimate with threshold level $t=C_{1}\frac{\log(p)}{n}$ for some constant $C_{1}$, then we have 
\begin{align}\label{eqn:DT:SIR}
\|P_{\widehat{\bbeta}_{DT}}-P_{\bbeta}\|^{2}\leq C_{2}\frac{s\log(p-s)}{n\lambda}
\end{align}
with probability at least $1-C_{4}\exp\left(-C_{3}\log(p)\right)$ for some positive constants $C_{2}, C_{3}$ and $C_{4}$.
\proof See the online supplementary file \citep{lin2016minimax}.
\end{thm}
From Theorem \ref{thm:DTSIR}, it is easy to see that, if $s=O(p^{1-\delta})$, the DT-SIR estimator $P_{\widehat{\bbeta}_{DT}}$ is rate optimal. 
Since there is a computational barrier for the rate optimal estimate of sparse PCA \citep{berthet2013computational},
the fact that the computationally efficient DT-SIR algorithm achieves the optimal rate suggests that sparse PCA might  not be  an appropriate prototype of SIR in high dimensions.

\section{Numerical Studies} In this section, we illustrate three aspects of the high dimensional behavior of SIR via numerical experiments. 
The first experiment focuses on the impacts of the choice of $H$ in SIR: the larger the $H$, the more accurate the estimate of eigenvalue of $var(\bbE[\vx|y])$. The second experiment aims at providing supporting evidence of Conjecture \ref{conj:derivative}. The third experiment demonstrates empirical performances of the  DT-SIR algorithm.

\subsection{Effects of  $H$}\label{subsec:numerical:H}
Our numerical results below show that the accuracy of estimating the  eigenvalues of $var(\bbE[\vx|y])$ depends on the choice of $H$.
Let us consider the following linear model:\footnote{Up to a monotone transform, this is the only case that we can give the explicit value of $\lambda(var(\bbE[\vx|y]))$.} 
\begin{align}
\mbox{Model $\mu$ :} \ & ~ y=\sqrt{\frac{\mu}{1-\mu}}\vx_{1}+\epsilon, \vx \sim N(0,\vI_{p}), \epsilon \sim N(0,1).
\end{align}
It is easy to see that the only non-zero eigenvalue of $var(\bbE[\vx|y])$ is $\mu$. 
The results are shown in Table~\ref{tab:sim2}, where  $H$ ranges in $\{2,5,10,50,100,200,500\}$, $\mu$  in $\{.5,.3,.1\}$ and $n$  in $\{5000,10000,50000,100000\}$. 
Each entry is the empirical mean (standard deviation), calculated based on 100 replications, of the SIR estimate of $\widehat{\mu}$ for given $\mu$ ,$n$ and $H$.

\begin{table}[H]
\begin{tabular}{|c|c|c|c|c|c|c|c|c|}
\hline
 & $n/1000$   & $H=2$ & $H=5$ & $H=10$  &  $H=50$ & $H=100$ & $H=200$ & $H=500$\\
\hline
\hline
\multirow{5}{*}{$\mu=.5$} 
& $5$  & 0.319 & 0.446 & 0.479 & 0.503 & 0.509 & 0.520 & 0.551\\ 
& & {\footnotesize(0.013)} & {\footnotesize(0.017)} & {\footnotesize(0.017)} & {\footnotesize(0.016)} & {\footnotesize(0.017)} & {\footnotesize(0.018)} & {\footnotesize(0.017)}
\\
& $10$   & 0.318 & 0.448 & 0.480 & 0.500 & 0.505 & 0.510 & 0.525 \\
& & {\footnotesize(0.009)} & {\footnotesize( 0.012)} & {\footnotesize( 0.012)} & {\footnotesize( 0.012)} & {\footnotesize( 0.012)} & {\footnotesize( 0.013)} & {\footnotesize( 0.012)}
\\ 
& $50$  & 0.319 & 0.448 & 0.479 & 0.498 & 0.500 & 0.501 & 0.504 \\
& &  {\footnotesize(0.004)} &  {\footnotesize(0.006)} &  {\footnotesize(0.005)} &  {\footnotesize(0.006)} &  {\footnotesize(0.005)} &  {\footnotesize(0.006)} &  {\footnotesize(0.006)}
\\
& $100$   & 0.319 & 0.448 & 0.479 & 0.498 & 0.499 & 0.501 & 0.503
\\
& & {\footnotesize(0.003)} & {\footnotesize( 0.004)} & {\footnotesize(0.004)} & {\footnotesize(0.004)} & {\footnotesize(0.004)} & {\footnotesize(0.004)} & {\footnotesize(0.004)}
\\[1pt]
\hline
\multirow{5}{*}{$\mu=.3$}  
& $5$  & 0.190 & 0.271 & 0.288 & 0.307 & 0.313 & 0.328 & 0.371 \\ 
& & {\footnotesize(0.011)} & {\footnotesize(0.012)} & {\footnotesize(0.014)} & {\footnotesize(0.015)} & {\footnotesize(0.015)} & {\footnotesize(0.014)} & {\footnotesize(0.016)}
\\
& $10$  & 0.191 & 0.27 & 0.288 & 0.302 & 0.307 & 0.312 & 0.335 \\
& & {\footnotesize(0.008)} & {\footnotesize(0.009)} & {\footnotesize(0.010)} & {\footnotesize(0.010)} & {\footnotesize(0.010)} & {\footnotesize(0.010)} & {\footnotesize(0.012)}
\\
& $50$    & 0.191 & 0.269 & 0.288 & 0.299 & 0.3 & 0.302 & 0.307\\
& & {\footnotesize(0.003)} & {\footnotesize(0.004)} & {\footnotesize(0.005)} & {\footnotesize(0.005)} & {\footnotesize(0.005)} & {\footnotesize(0.005)} & {\footnotesize( 0.004)}
\\
& $100$  & 0.191 & 0.269 & 0.288 & 0.299 & 0.3 & 0.301 & 0.303 \\
& & {\footnotesize(0.002)} & {\footnotesize(0.003)} & {\footnotesize(0.003)} & {\footnotesize( 0.004)} & {\footnotesize( 0.003)} & {\footnotesize( 0.003)} & {\footnotesize( 0.003)}
\\[1pt]
\hline
\multirow{5}{*}{$\mu=.1$}  
& $5$  
& 0.064 & 0.091 & 0.098 & 0.109 & 0.117 & 0.136 & 0.190 \\ 
& & {\footnotesize(0.007)} & {\footnotesize(0.008)} & {\footnotesize(0.009)} & {\footnotesize(0.009)} & {\footnotesize(0.009)} & {\footnotesize(0.010)} & {\footnotesize(0.010)}
\\
& $10$   & 0.0643 & 0.0901 & 0.0973 & 0.103 & 0.108 & 0.117 & 0.144 \\
& & {\footnotesize(0.005)} & {\footnotesize(0.006)} & {\footnotesize(0.006)} & {\footnotesize(0.006)} & {\footnotesize(0.006)} & {\footnotesize(0.006 )}& {\footnotesize(0.007)}
\\
& $50$   & 0.0638 & 0.0899 & 0.0963 & 0.101 & 0.101 & 0.103 & 0.109
\\
& & {\footnotesize(0.002)} & {\footnotesize( 0.003)} & {\footnotesize( 0.003)} & {\footnotesize( 0.003)} & {\footnotesize( 0.003)} & {\footnotesize( 0.003)} & {\footnotesize( 0.003)}
\\
& $100$  & 0.0636 & 0.0898 & 0.0961 & 0.100 & 0.100 & 0.102 & 0.104 \\
& & {\footnotesize( 0.001)} & {\footnotesize( 0.002)} & {\footnotesize( 0.002)} & {\footnotesize( 0.002)} & {\footnotesize( 0.002)} & {\footnotesize( 0.002)} & {\footnotesize( 0.002)}
\\
\hline
\end{tabular}
\caption{The empirical mean (standard error) of the SIR estimate   $\widehat{\lambda}(\mu)$ for $\mu$}\label{tab:sim2}
\end{table}
From Table \ref{tab:sim2}, it is clear that the larger the $H$ is, the more accurate estimation of the eigenvalue is. 
Cautious reader may notice that, in the row with $\mu=.1$ and $n=5000$, the empirical mean and the standard error are not behaving as we have expected, e.g., when $H=500$, the empirical mean and standard error are 0.190 and 0.010, respectively, which are worse than the case with $H=10$ (or $50$). 
This is not contradicting our theory. Note that in the Lemma \ref{lem:main:deviation}, the deviation property of $\widehat{\lambda}$ depends on the value $\frac{n\mu}{H^{2}}$, i.e., the larger the $\frac{n\mu}{H^{2}}$ is, the more concentrated the $\widehat{\lambda}$ is. In particular, for the entry corresponding to $\mu=0.1$, $n/1000=5$ and $H=500$,  the value $\frac{n\mu}{H^{2}}=1/500$ is much smaller than the corresponding value, 5, associated with the entry with $\mu=0.1$, $n/1000=5$ and $H=10$.

\subsection{Support Evidences of Conjecture \ref{conj:derivative}}\label{subsubsection:numerical:conjecture}
%In this subsection, we provide some support evidences of the Conjecture \ref{conj:derivative}.
Let us consider the following model with two indexes:
\begin{align}\label{model:simulation:conj}
\mbox{Model $\mu$ :}& ~ y=
\sqrt{\mu}(1+g(\vx_{1}))(g(\vx_{1})+g(\vx_{2}))+\epsilon
\end{align}
where $g:\bbR\mapsto \bbR$ is a smooth function such that for a small constant $\delta>0$,
\begin{align}
g(x)=\begin{cases}
x & \mbox{ if } |x|\leq 100-\delta\\
0 & \mbox{ if } |x| \geq 100+\delta
\end{cases}
\end{align}
and $|g'(x)|\leq C$ for some constant $C$.
Let $\lambda_{1}(\mu)$ and $\lambda_{2}(\mu)$ be the two eigenvalues of $var(\bbE[\vx|y])$. 
Since we know that the absolute value of the derivative of the link function $\leq C\sqrt{\mu}$, we want to check if $C_{1}\mu \leq \lambda_{2}(\mu) \leq \lambda_{1}(\mu) \leq C_{2}\mu$ holds for some positive constant $C_{1}$ and $C_{2}$ and if model \eqref{model:simulation:conj} belongs to $\mathcal{F}_{2}(C_{1}\mu, C_{2}/C_{1})$.
We study the boundedness of $\lambda_{1}(\mu)/\mu$ and $\lambda_{2}(\mu)/\mu$ via numerical simulation. In the simulation, we choose $H$ to be 20. Let $\mu$ range in $\{1,.5,.1,.05,.01,.005,.001\}$ and $n$ range in $\{10^{3},10^{5},10^{5},10^{6}\}$.
\begin{table}[H]
\begin{tabular}{|c|c|c|c|c|c|c|c|c|}
\hline
 & $n$   & $\mu=1$ & $\mu=.5$ & $\mu=.1$  &  $\mu=.05$ & $\mu=.01$ & $\mu=.005$ & $\mu=.0001$\\
\hline
\hline
\multirow{4}{*}{$\lambda_{1}(\mu)/\mu$}  & $n=10^{3}$  & 0.3358 &    0.5969   &  1.6333   &  2.1297 &    4.5681   &  7.1434  &  30.6206\\ 
& $n=10^{4}$   & 0.3276 &    0.5676   &  1.4416 &    1.7627  &   2.1511   &  2.3908  &   4.3137 \\
& $n=10^{5}$   &  0.3272  &  0.5650   & 1.4153  &  1.7092  &  1.9662 &   2.0006   & 2.2554 \\
& $n=10^{6}$   &   0.3268  &  0.5651   & 1.4125 &   1.7052   & 1.9465   & 1.9695  &  1.9780 \\
\hline
\multirow{3}{*}{$\lambda_{2}(\mu)/\mu$}  & $n=10^{3}$  &  0.1068 &    0.1436   &  0.3227 &    0.5381 &    2.3159   &  4.2023 &   19.5854 \\ 
& $n=10^{4}$   &  0.0899 &   0.1087   & 0.1248  &  0.1206   & 0.2701 &   0.5061    & 2.3384 \\
& $n=10^{5}$   &  0.0899    &  0.1059    &  0.1014   &   0.0840  &    0.0462    &  0.0620   &   0.2366
\\
& $n=10^{6}$   &    0.0898 &     0.1063    &   0.1001   &   0.0795  &    0.0297   &   0.0190  &    0.0278 \\
\hline
\end{tabular}
\caption{The empirical expectation of  $\lambda_{i}(\mu)/\mu$ }\label{tab:sim1}
\end{table}
In Table \ref{tab:sim1}, each entry is the average of 100 replications. For fixed $\mu$, the larger $n$, the more accurate estimation of $\lambda_{i}(\mu)/\mu, i=1,2.$ In particular, it is easy to see from the row with $n=10^{6}$ that $\lambda_{i}(\mu)/\mu, i=1,2.$ are bounded. 
The row with $n=10^{3}$ seems to be contradicting to our conjecture \ref{conj:derivative}, where $\lambda_{i}(\mu)/\mu, i=1,2,$ might be diverging as $\mu \rightarrow 0$.  This is actually not a contradiction, since  we know that the deviation property of $\lambda(var(\bbE[\vx|y]))$ depends on the product $n\lambda$ from Lemma \ref{lem:main:deviation}. Thus, to get accurate estimate of $\lambda(var(\bbE[\vx|y]))$, we require more samples if $\lambda$ is small.

\subsection{Performance of DT-SIR} 
\begin{figure}[H]
  \centering
  \caption{Average values of $\|P_{\widehat{\bbeta}_{DT}}-P_{\bbeta}\|^{2}$}\label{fig:loss}
    \includegraphics[width=0.7\textwidth]{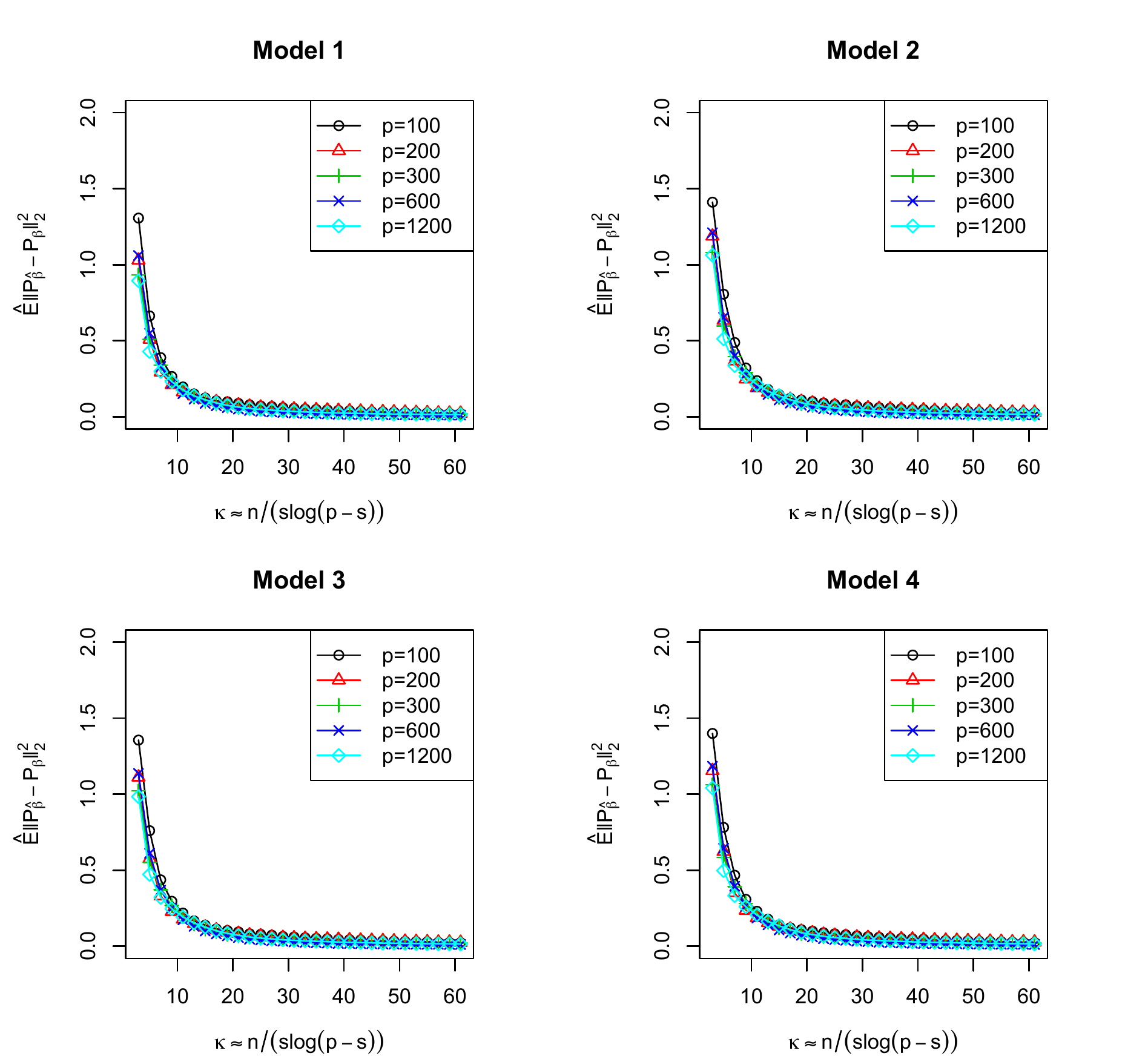}      
\end{figure}

In this section, we assume the exact sparsity $s=O(p^{1-\delta})$ for some $\delta\in (0,1)$, and consider the following data generating models,
\begin{align*}
\text{Model 1}: y =& \boldsymbol{x}^\tau\boldsymbol{\beta}+\sin(\boldsymbol{x}^\tau \boldsymbol{\beta}) + \epsilon,\\
\text{Model 2}: y = & 2\arctan(\boldsymbol{x}^\tau\boldsymbol{\beta})+\epsilon,\\
\text{Model 3}: y = & (\boldsymbol{x}^\tau\boldsymbol{\beta})^3+\epsilon,\\
\text{Model 4}: y = & \sinh (\boldsymbol{x}^\tau\boldsymbol{\beta})+\epsilon,
\end{align*}
where $\boldsymbol{x}\sim N(\boldsymbol{0}, \boldsymbol{I}_p)$, $\epsilon \sim N(0,1)$, $\boldsymbol{x} \independent \epsilon$, and
$\bbeta$ is a fixed vector with $s$ nonzero coordinates. 
Let $\kappa = \{s\log(p-s)/n\}^{-1}$.
The dimension $p$  of the predictors takes value in $\{100, 200,$ $300, 600, 1200\}$, the sparsity parameter $\delta$ is fixed at $0.5$, and $\kappa$ takes values in $\{3,5,7,\ldots, 61\}$. For each $(p,\kappa)$ combination, $s=\lfloor p^{1-\delta}\rfloor$, $n =\lfloor \kappa s \log(p-s) \rfloor$, and we simulate data from each model 1000 times. We then get the estimate $\widehat{\bbeta}_{DT}$ using DT-SIR algorithm, and the results of the average values of $\|P_{\widehat{\bbeta}_{DT}}-P_{\bbeta}\|^{2}$ for each model with each  $(p,\kappa)$ combination are shown in Figure \ref{fig:loss}, which shows the distance between the estimated projection matrix and the true one becomes smaller as $\kappa$ increases for all fixed $p$.

\begin{figure}
  \centering
  \caption{Average values of $\kappa \|P_{\widehat{\bbeta}_{DT}}-P_{\bbeta}\|^{2}$}\label{fig:kappa_loss}
    \includegraphics[width=0.7\textwidth]{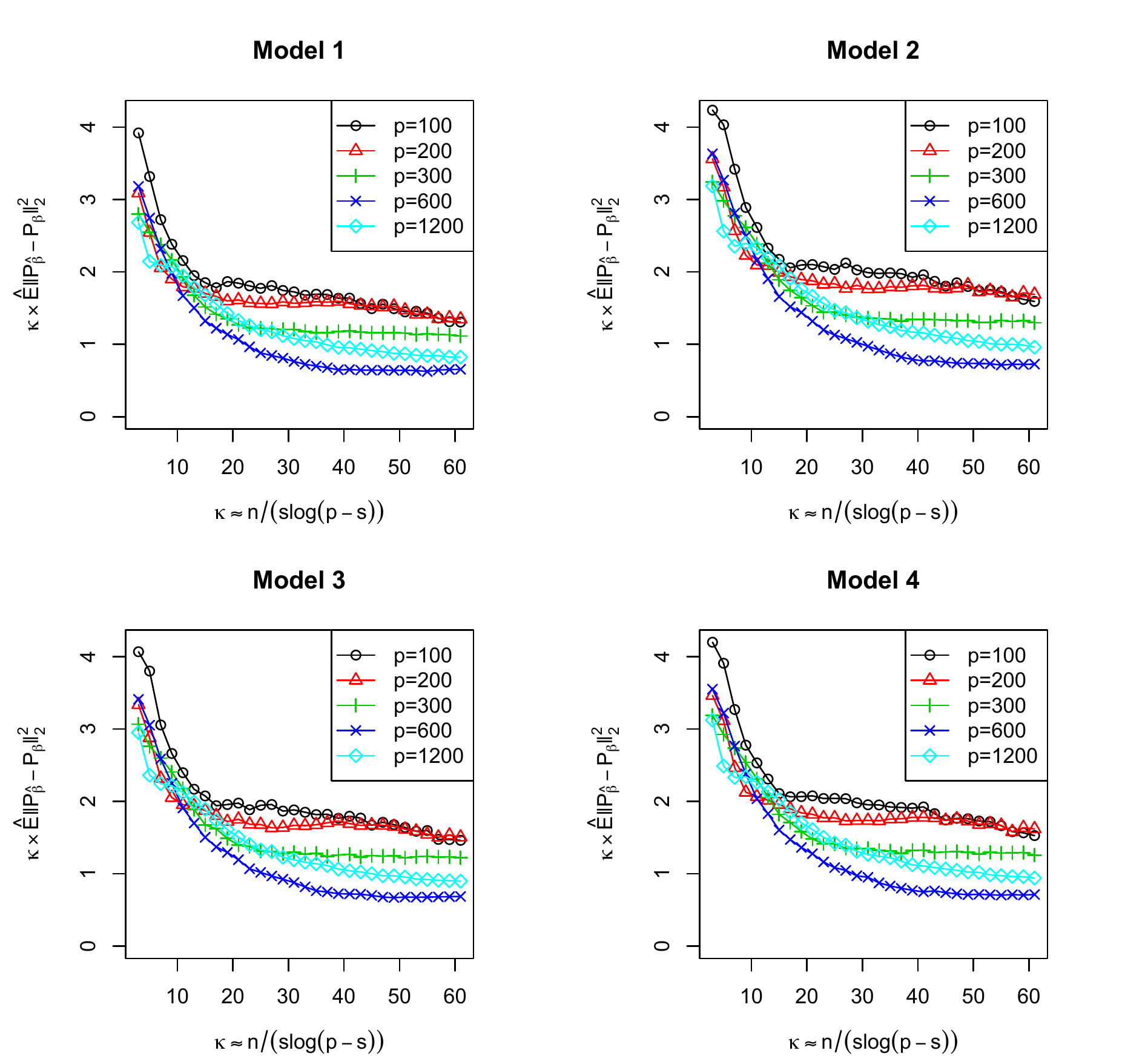}      
\end{figure}

According to Theorem \ref{thm:DTSIR}, $\kappa \|P_{\widehat{\bbeta}_{DT}}-P_{\bbeta}\|^{2}$ is less than a constant with high probability. Therefore, we also  the average values of $\kappa * \|P_{\widehat{\bbeta}_{DT}}-P_{\bbeta}\|^{2}$ for these models in Figure \ref{fig:kappa_loss}, which demonstrates that $\kappa \|P_{\widehat{\bbeta}_{DT}}-P_{\bbeta}\|^{2}$ is a decreasing function of $\kappa$ and tends to be stable when $\kappa$ becomes large enough. These empirical results also validate Theorem \ref{thm:DTSIR}.

\newpage

\section{Proofs}
We need the following technical lemma, which is a direct corollary of the {\it `key lemma' } in \cite{lin2015consistency} :
 
\begin{lem} \label{lem:main:deviation}
Assume that $f \in \mathcal{F}_{d}(\lambda,\kappa)$ in the model \eqref{model:modified:multiple}. Let $\widehat{\bLambda}_{H}$ be the SIR estimate \eqref{eqn:lambda} of $var(\bbE[\vx|y])(=\bLambda)$.
There exist positive absolute constants $C$, $C_{1},C_{2}$ and $C_{3}$ such that, for any $f \in \mathcal{F}_{d}(\lambda,\kappa)$ and any $\nu >1$, if  $H>C(\nu^{1/\vartheta}\vee d)$ for sufficiently large constant $C$, then for any unit vector $\bbeta$ that lies in  the column space of $\bLambda$, we have
\begin{align}\label{eqn:inline:temp}
  \Big| \bbeta^{\tau}\left(\widehat{\bLambda}_{H}-\bLambda\right)\bbeta\Big| > \frac{1}{2\nu}\bbeta^{\tau}\bLambda\bbeta 
\end{align}  
with probability at most
\begin{align*}
 C_{1}\exp\left(-C_{2}\frac{n\bbeta^{\tau}\bLambda\bbeta}{H^{2}\nu^{2}} +C_{3}\log(H)\right).
\end{align*}
In particularly, if $d$ and $\nu$ are bounded, we can choose $H$ to be a large enough  finite integer such that \eqref{eqn:inline:temp} holds with high probability.
\proof It is a direct corollary of the {\it `key lamma'} in \cite{lin2015consistency}. \epf
\end{lem} 

\subsection{Proof of Theorem \ref{thm:risk:oracle:upper:d} }\label{sec:oracle:risk:proof}

Suppose that we have  $n=Hc$ samples $(y_{i},\vx_{i})$ from the distribution defined by the model  $\mathcal{M}=(\vV,f) \in \mathfrak{M}(p,d,\kappa,\lambda)$. 
Let $H=H_{1}d$ where $H_{1}$ is a  sufficiently large integer and $\widehat{\vV}=(\widehat{\vV}_{1},...,\widehat{\vV}_{d})$ where $\widehat{\vV}_{i}$ is the eigen-vector associated to the $i$-th largest eigen-value of $\widehat{\bLambda}_{H}$.  
We introduce the following decomposition 
\begin{equation*}
\begin{aligned}\label{decomposition:all}
\vx&=P_{\mathcal{S}}\vx+P_{\mathcal{S}^{\perp}}\vx
%=\bbE[\vx|y]+(P_{\mathcal{S}}\vx-\bbE[\vx|y])+P_{\mathcal{S}^{\perp}}\vx\\
\triangleq \vz+\vw,
\end{aligned}
\end{equation*}
i.e., $\vz$ lies in the central space $\mathcal{S}$ and $\vw$ lies in the space $\mathcal{S}^{\perp}$ which is perpendicular to $\mathcal{S}$. %Sometimes, we will use $\vz\triangleq\vm+\vv$ and $\bdelta\triangleq\vv+\vw$ to denote $P_{\mathcal{S}}\vx$ and the noise does not lie in the central curve respectively.  
Let $\vV^{\perp}$ be a $p\times (p-d)$ orthogonal matrix such that $\vV^{\tau}\vV^{\perp}=0$. 
Since $\mathcal{S}=span\{\vV\}$ and $\vx\sim N(0,\vI_{p})$, we may write $\vw=\vV^{\perp}\bepsilon$ for some $\bepsilon \sim N(0,\vI_{p-d})$.  Thus we know that $\bSigma_{\vw}\triangleq var(\vw)=\vV^{\perp}\vV^{\perp,\tau}$.
We introduce
the notation $\overline{\vz}_{h,\cdot}$, $\overline{\vw}_{h,\cdot}$, 
%$\overline{\bdelta}_{h,\cdot}$, $\overline{\overline{\bdelta}}$, 
and $\overline{\bepsilon}_{h,\cdot}$ similar to
the definition of $\overline{\vx}_{h,\cdot}$.
Let
%\begin{equation}\label{eqn:notation:keynotation}
%\begin{aligned}
$\bold{\mathcal{Z}}=\frac{1}{\sqrt{H}}
\left(\overline{\vz}_{1,\cdot} ~,~ \overline{\vz}_{2,\cdot},...,~ 
\overline{\vz}_{H,\cdot}\right)$, 
$\bold{\mathcal{W}}=\frac{1}{\sqrt{H}}
\left(\overline{\vw}_{1,\cdot} ~,~ \overline{\vw}_{2,\cdot},...,~ \overline{\vw}_{H,\cdot}\right)$,
$\bold{\mathcal{E}}=\frac{1}{\sqrt{H}}
\left(\overline{\bepsilon}_{1,\cdot} ~,~ \overline{\bepsilon}_{2,\cdot},...,~ \overline{\bepsilon}_{H,\cdot}\right)$
%\end{aligned}
%\end{equation}
be three $p\times H$ matrices formed by the vectors $\frac{1}{\sqrt{H}}\overline{\vz}_{h,\cdot}$, $\frac{1}{\sqrt{H}}\overline{\vw}_{h,\cdot}$,  and $\frac{1}{\sqrt{H}}\overline{\bepsilon}_{h,\cdot}$.
We have the following decomposition
\begin{equation}
\begin{aligned} \label{eqn:estimator:decomposition}
\widehat{\bLambda}_{H}&=\bold{\mathcal{Z}}\bold{\mathcal{Z}}^{\tau}
+\bold{\mathcal{Z}}\bold{\mathcal{W}}^{\tau}
+\bold{\mathcal{W}}\bold{\mathcal{Z}}^{\tau}
+\bold{\mathcal{W}}\bold{\mathcal{W}}^{\tau}\\
%&=\bLambda_{u}+\bold{\mathcal{Z}}\bold{\mathcal{W}}^{\tau}
%+\bold{\mathcal{W}}\bold{\mathcal{Z}}^{\tau}
%+\bold{\mathcal{W}}\bold{\mathcal{W}}^{\tau}\\
&=\bLambda_{u}+\bold{\mathcal{Z}}\bold{\mathcal{E}}^{\tau}\vV^{\perp,\tau}
+\vV^{\perp}\bold{\mathcal{E}}\bold{\mathcal{Z}}^{\tau}
+\vV^{\perp}\bold{\mathcal{E}}\bold{\mathcal{E}}^{\tau}\vV^{\perp,\tau}
\end{aligned}
\end{equation}
where we define $\bLambda_{u}\triangleq \mathcal{Z}\mathcal{Z}^{\tau}$ and use the fact $\mathcal{W}=\vV^{\perp}\mathcal{E}$. 
Since $\bepsilon \sim N(0,\vI_{p-d})$, we know that the entries $\mathcal{E}_{i,j}$ of $\mathcal{E}$ are $i.i.d.$ samples of $N(0,\frac{1}{n})$. 
First, we have the following lemma.
\begin{lem} \label{lem:minimx:key}
Let $\rho=\frac{p}{n}$. Assume that $\frac{p}{n\lambda}$ is sufficiently small. We have the following statements. 
%Let $\alpha=\rho\vee \frac{\log(n\lambda)}{n}  $
\begin{itemize}
\item[i)] There exist constants $C_{1}, C_{2}$ and $C_{3}$ such that
\end{itemize}
\begin{align*}
\bbP(\|\mathcal{W}\mathcal{W}^{\tau}\|> C_{1}\rho)\leq C_{2}\exp\left(-C_{3}p\right).
\end{align*}
\begin{itemize}
\item[ii)] For any vector $\bbeta\in \mathbb{R}^{p}$ and any $\nu>1$, 
let $E_{\bbeta}(\nu)=\Big\{~ \Big| \bbeta^{\tau}\left(\bLambda_{\vu}-\bLambda\right)\bbeta\Big| > \frac{1}{2\nu}\bbeta^{\tau}\bLambda\bbeta\Big\}$. 
Recall that $H=dH_{1}$. If we choose $H_{1}$ sufficiently large such that $H^{\vartheta}>C\nu$ for some positive constant $C$, there exist positive constants $C_{1}$, ..., $C_{3}$ and $C_{4}$ such that 
\end{itemize}
\begin{align*}
\bbP\left(  \bigcup_{\bbeta}E_{\bbeta} (\nu)\right)\leq C_{1}\exp\left(-C_{2}\frac{n\lambda}{H^{2}\nu^{2}} +C_{3}\log(H)+C_{4}d\right).
\end{align*} 
\begin{itemize}
\item[iii)]  For any $\nu>1$, there exist positive constants $C_{1}$,..., $C_{6} $ and  $C_{7}$, such that 
\end{itemize}
\begin{align*}
\bbP\left(\|\mathcal{W}\mathcal{Z}^{\tau}\|> C_{7}\sqrt{\kappa\lambda\rho} \right)\leq& 
C_{1}\exp\left(-C_{2}\frac{n\lambda}{H^{2}\nu^{2}} +C_{3}\log(H)+C_{4}d\right)\\
&+C_{5}\exp\left(-C_{6}p\right).
\end{align*}

\proof 
$i)$ is a direct corolllary of Lemma \ref{random:nonasymptotic}.  %In Lemma \ref{random:nonasymptotic}, by choosing $t^{2}=p$ , we have
%\begin{align}
%\bbP(\|\mathcal{W}\mathcal{W}^{\tau}\|>9\rho) \leq \exp\left(-p \right).
%\end{align}
$ii)$ is a direct corollary of Lemma \ref{lem:main:deviation} and the usual $\epsilon$-net argument.
$iii)$ is a direct corollary of $i)$ and $ii)$
\epf
\end{lem}

Let $\ttE=\ttE_{1}\cap \ttE_{2} \cap \ttE_{3}$ where $\ttE_{1}=\Big\{~ \|\bold{\mathcal{W}}\bold{\mathcal{W}}^{\tau}\|\leq C\rho~\Big\}$, 
$\ttE_{2}=\Big\{~ \|\bold{\mathcal{W}}\bold{\mathcal{Z}}^{\tau}\|\leq 4 \sqrt{\kappa\lambda\rho} ~\Big\}$, 
%$\ttE_{3}=\Big\{~ \|\bLambda_{u}-\bLambda\|\leq \frac{1}{2\nu}\|\bLambda\| ~\Big\}$ and $\ttE=\ttE_{1}\cap \ttE_{2} \cap \ttE_{3}$.  
$\ttE_{3}=\Big\{~ \|\bLambda_{u}-\bLambda\|\leq \frac{1}{2\nu}\kappa\lambda~\Big\}$.

\begin{cor}\label{cor:minimax:key} Lemma \ref{lem:minimx:key} implies the following simple results where $C$ stands for some absolute constant which might be varying in different statements.  
\begin{itemize}
\item[$a)$]
If $n\lambda \leq e^{p}$, we have $\bbP\left( \ttE^{c}\right) \leq \frac{CH^{2}}{n\lambda}$.
\item[$b)$]Conditioning on $\ttE_{3}$,  we have  $\lambda_{d}(\bLambda_{\vu}) \geq (1-\frac{\kappa}{2\nu})\lambda $.
\item[$c)$] Conditioning on $\ttE$, if $\frac{p}{n\lambda}$ is sufficiently small, we have 
$\|\widehat{\bLambda}_{H}-\bLambda_{\vu}\|
%&\leq \| %\bold{\mathcal{W}}\bold{\mathcal{W}}^{\tau} %\|
%+\| %\bold{\mathcal{W}}\bold{\mathcal{Z}}^{\tau} %\|
%+\| %\bold{\mathcal{Z}}\bold{\mathcal{W}}^{\tau} %\| \\
%&\leq 9\rho + 8\sqrt{\kappa \lambda \rho}
 \leq C\sqrt{\frac{\kappa \lambda p}{n} }$. 
\item[$d)$] Conditioning on $\ttE$, If $\frac{p}{n\lambda}$ is sufficiently small, we have
$\lambda_{d+1}(\widehat{\bLambda}_{H}) < \frac{1}{4}\lambda$.
\end{itemize}
\end{cor}

Now we start the proof of Theorem \ref{thm:risk:oracle:upper:d}. Note that
\begin{equation*}
\begin{aligned}
\bbE\|\widehat{\vV}\widehat{\vV}^{\tau}&-\vV\vV^{\tau}\|^{2}_{F}\\
&=\underbrace{\bbE\|\widehat{\vV}\widehat{\vV}^{\tau}-\vV\vV^{\tau}\|^{2}_{F}\bold{1}_{\ttE}}_{I}
+
\underbrace{\bbE\|\widehat{\vV}\widehat{\vV}^{\tau}-\vV\vV^{\tau}\|^{2}_{F}\bold{1}_{\ttE^{c}}}_{II}.
\end{aligned}
\end{equation*}
\paragraph{For $II$}
It is easy to see that 
\begin{align*}
II \leq 2(d\wedge (p-d))\bbP(\ttE^{c})=2d\bbP\left( \ttE^{c}\right) \leq \frac{CdH^{2}}{n\lambda}=\frac{Cd^{3}H_{1}^{2}}{n\lambda}.
\end{align*}
\paragraph{For $I$}  
 Let 
 $ \bLambda_{\vu}=
  \widetilde{\vV}\vD_{H}\widetilde{\vV}^{\tau}$
be the spectral decomposition of $\bLambda_{\vu}$, where $\widetilde{\vV}$ is a $p\times d$ orthogonal matrix and $\vD_{H}$ is a $d\times d$ diagonal matrix.
Conditioning on $\ttE$, we know that  $\widetilde{\vV}$ and $\vV$ are sharing the same column space. 
Thus we have  $\widetilde{\vV}\widetilde{\vV}^{\tau}=\vV\vV^{\tau}$. 
%Recall that $\widehat{\vV}_{H}$ and $\vV$ are the top eigen-space of $\widehat{\bLambda}_{H}$ and $\bLambda_{\vu}$ respectively. 
Let us apply  the Sin-Theta theorem (e.g., Lemma \ref{lem:sin_theta}) to  the pair of symmetric matrices
$(\bLambda_{\vu}, \widehat{\bLambda}_{H}=\bLambda_{\vu}+\vQ$) where $\vQ\triangleq\widehat{\bLambda}_{H}-\bLambda_{u}$.
Since $\frac{p}{n\lambda}$ is sufficiently small, conditioning on $\ttE$, we have $\lambda_{d+1}(\widehat{\bLambda}_{H})\leq \frac{1}{4}\lambda$ and $\lambda_{d}(\bLambda_{\vu})=\lambda_{d}(\vD_{H})\geq \frac{\lambda}{2}$. Thus, 
%\begin{align}
%\bLambda_{\vu}&=[\widetilde{\vV},%\widetilde{\vV}^{\perp}]\left[ \begin{array}{cc}
%\vD_{H} & 0\\
%0 & 0
%\end{array} \right]
%\left[ \begin{array}{c}
%\widetilde{\vV}^{\tau}\\
%\widetilde{\vV}^{\perp,\tau}
%\end{array} \right]
%\end{align}
%and
%\begin{equation}
%\begin{aligned}
%\widehat{\bLambda}_{H}&=\bLambda_{\vu}+\left(\widehat{\bLambda}_{H}-
%\bLambda_{\vu}\right)\\
%&=[\widehat{\vV},\widehat{\vV}^{\perp}]\left[ %\begin{array}{cc}
%\vD_{1} & 0\\
%0 & \vD_{2}
%\end{array} \right]
%\left[ \begin{array}{c}
%\widehat{\vV}^{\tau}\\
%\widehat{\vV}^{\perp,\tau}
%\end{array} \right]
%\end{aligned},
%\end{equation}
we have
\begin{align*}
&\bbE\|\vV\vV^{\tau}-\widehat{\vV}\widehat{\vV}^{\tau}\|^{2}_{F}\bold{1}_{\ttE}=\bbE\|\widetilde{\vV}\widetilde{\vV}^{\tau}-\widehat{\vV}\widehat{\vV}^{\tau}\|^{2}_{F}\bold{1}_{\ttE}\\
\leq &\frac{32}{\lambda^{2}}
\min
\left(
\bbE\|\widetilde{\vV}^{\perp,\tau} \vQ \widehat{\vV}\|_{F}^{2}\bold{1}_{\ttE}, \bbE\|\widetilde{\vV}^{\tau}\vQ \widehat{\vV}^{\perp}\|_{F}^{2}\bold{1}_{\ttE}
\right)\\
 \leq &\frac{32}{\lambda^{2}}\min\left(\bbE\|\vQ\widetilde{\vV}\|_{F}^{2}\bold{1}_{\ttE}, \bbE\|\vQ\widetilde{\vV}^{\perp}\|_{F}^{2}\bold{1}_{\ttE}\right).
\end{align*}
Since $\widetilde{\vV}$ and $\vV$ are sharing the same column space, we have
$
\widetilde{\vV}^{\tau}\mathcal{W}=\vV^{\tau}\mathcal{W}=0$ and  $\widetilde{\vV}^{\perp,\tau}\mathcal{Z}=\vV^{\perp,\tau}\mathcal{Z}=0$. 
Thus, we have
\begin{align*}
\widetilde{\vV}^{\tau}\vQ
&=\widetilde{\vV}^{\tau}\mathcal{Z}\mathcal{W}^{\tau},\quad
\widetilde{\vV}^{\perp,\tau}\vQ
=\widetilde{\vV}^{\perp,\tau}\mathcal{W}\mathcal{W}^{\tau}
+\widetilde{\vV}^{\perp,\tau}\mathcal{W}\mathcal{Z}^{\tau}.
\end{align*}
Conditioning on $\ttE$, we have $\|\bLambda_{u}\|_{2}\leq 2\kappa\lambda$. 
Thus
\begin{align*}
&\min\left(\bbE\| \vQ\widetilde{\vV}\|_{F}^{2}\bold{1}_{\ttE}, \bbE\|\vQ\widetilde{\vV}^{\perp}\|_{F}^{2}\bold{1}_{\ttE}\right)
\leq 2\bbE\|\widetilde{\vV}^{\tau}\mathcal{Z}\mathcal{W}^{\tau}\|^{2}_{F}\bold{1}_{\ttE}\leq 4\kappa\lambda \bbE\|\mathcal{W}^{\tau}\|^{2}_{F}\leq \frac{4\kappa\lambda}{n} d(p-d).
\end{align*}
Since $\kappa$ is assumed to be fixed, we know that if $\frac{p}{n\lambda}$ is sufficiently small and $d^{2}\leq p$, we have
\begin{align*}
\sup_{\mathcal{M} \in \mathfrak{M}(p,d,\kappa,\lambda)}\bbE\|\widehat{\vV}\widehat{\vV}^{\tau}-\vV\vV^{\tau}\|^{2}_{F}\prec \frac{d(p-d)}{n\lambda}. 
\end{align*}
\epf

\subsection{Proof of the Theorem \ref{thm:risk:sparse:upper:d}} \label{sec:sparse:risk:proof}

Before we start proving this Theorem, we need some preparations. First, the following lemmas will be used frequently during the proofs. 
\begin{lem}\label{inline:trivial1}
Let $\vK$ be an $a\times b$ matrix with each entry being i.i.d. standard normal random variables. Then, we have
$
\bbE[\|\vK\vK^{\tau}\|^{2}_{F}]=ab(a+b+1)$
 and $
\bbE[\|\vK\|^{2}_{F}]=ab
$.
\proof It follows from elementary calculations.
\epf
\end{lem}
\begin{lem}\label{lem:elementary:trivial2}
Let $A$ , $B$ be  $l\times m$ and  $m\times n$ matrices, respectively, we have
$
\|AB\|_{F} \leq \|A\|_{2}\|B\|_{F},
$ where $\|A\|_{2}$ denotes the largest singular value of $A$. 
\proof It follows from elementary calculations.
\epf
\end{lem}

\begin{lem}\label{lem:elementary:trivial3}
Let $A$, $B$ be $m\times l$ orthogonal matrices, i.e., $A^{\tau}A=I_{l}=B^{\tau}B$, and let $M$ be an $l\times l$ positive definite matrix with eigenvalues $d_j$ such as $0<\lambda\leq d_{l} \leq d_{l-1} \leq ... \leq d_{1} \leq \kappa \lambda$. If $A^{\tau}B$ is a diagonal matrix with non-negative entries % and $\|A^{\tau}B-I_{l}\|_{F}$ is sufficiently small
, then there exists a constant $C$ which only depends on $\kappa$ such that 
$
\|AMA^{\tau}-BMB^{\tau}\|_{F}\leq C\lambda \|AA^{\tau}-BB^{\tau}\|_{F}.$
\proof 
Let $\Delta=I_{l}-B^{\tau}A$, then $0 \leq \Delta_{ii}\leq 1$ for $1\leq i \leq l$.  If $C>2\kappa^{2}-1$,  we have
\begin{equation}\nonumber
\begin{aligned}
\|AMA^{\tau}-BMB^{\tau}\|_{F}^{2}&=2tr(M^{2}\Delta)-tr(M\Delta M\Delta)\leq 2\kappa^{2}\lambda^{2} tr(\Delta)-\lambda^{2}tr(\Delta^{2})\\ 
&\leq C\lambda^{2}(2tr(\Delta)-tr(\Delta^{2}))=C\lambda^{2}\|AA^{\tau}-BB^{\tau}\|_{F}^{2}.
\end{aligned}
\end{equation}
 \epf

\end{lem}

\begin{lem}\label{cor:elemetary:temp}
For a positive definite matrix M with eigenvalue  $\lambda_{1}\geq ... \geq \lambda_{d}>0$ and orthogonal matrices A,B,E,F, i.e., $A^{\tau}A=B^{\tau}B=E^{\tau}E=F^{\tau}F=I_{d}$, 
we have
\[
\frac{\lambda_{d}}{2}\|AB^{\tau}-EF^{\tau}\|^{2}_{F} \leq ~\langle AMB^{\tau},AB^{\tau}-EF^{\tau} \rangle ~\leq \frac{\lambda_{1}}{2}\|AB^{\tau}-EF^{\tau}\|^{2}_{F}.
\]
\proof It is a direct corollary of the Lemma 8 in \cite{gao2014minimax}. \epf
\end{lem}

\begin{lem}[Sparse approximation]  \label{lem:risk:sparse:approximation}
Let $\vV \in \mathbb{O}_{s,q}(p,d)$ and $k \in [p]$, where $\mathbb{O}_{s,q}(p,d)$ is defined near \eqref{def:sparse orthogonal}.
 Let $\|\vV_{(i)*}\|$ denote its i-th largest row norm. Then
\begin{align}
\sum_{i>k}\|\vV_{(i)*}\|^{2} &\leq \frac{q}{2-q}k(s/k)^{2/q}.
\end{align}
In particular, if $k$ is chosen to be $k_{s,q}$ defined near \eqref{def:sparse orthogonal}, we know that 
\begin{align}
\sum_{i>k}\|\vV_{(i)*}\|^{2} &\leq \frac{q}{2-q}\epsilon^{2}_{n}.
\end{align}
\proof This is a direct corollary of the Lemma 7 in \cite{cai2013sparse}.\epf
\end{lem}

%\begin{lem}[Sparse approximation]  \label{lem:risk:sparse:approximation}
%Let $\vV \in \mathcal{F}_{q}(s,p)$ and $k \in [p]$, where $\mathcal{F}_{q}(s,p)$ defined . Let $\|\vV_{(i)*}\|$ denote its i-th largest row norm. Then
%\begin{align}
%\sum_{i>k}\|\vV_{(i)*}\|^{2} &\leq \frac{q}{2-q}k(s/k)^{2/q}\\
%\|\vV_{S}\vV_{S}^{\tau}-\vV\vV^{\tau}\|^{2}_{F} &\leq \frac{Cq}{2-q}\epsilon^{2}_{n}\\
%\|\vV_{S}\vD\vV_{S}^{\tau}-\vV\vD\vV^{\tau}\|^{2}_{F} &\leq \frac{Cq}{2-%q}\lambda^{2}\epsilon^{2}_{n}
%\end{align}
%where $\epsilon_{n}^{2} \triangleq \frac{1}{n\lambda}\left(dk+k\log\frac{ep}{k} \right)$ is assumed to be sufficiently small.
%\end{lem}

\begin{lem} \label{sparse:random}
Let $\bSigma=\vV\vD\vV^{\tau}$ be a $p\times p$ positive semidefinite matrix where $\vV$ is a $p\times d$ orthogonal matrix and $\vD$ is a $d\times d$ diagonal matrix with entries 
$\lambda \leq d_{d}\leq ... \leq d_{1} \leq \kappa \lambda$.
For a subset $S$ of indices with $|S|=k$, 
let $J_{S}$ be a diagonal matrix such that $J_{S}(i,i)=1$  if $i \in S$ and $J_{S}(i,i)=0$ if $i \not \in S$ . Let $\bSigma_{S}=J_{S}\bSigma J_{S}$ and let $\bSigma_{S}=\vV_{1}D_{1}\vV_{1}^{\tau}$ be the eigen-decomposition of $\bSigma_{S}$. We have
\begin{align*}
\|\bSigma-\bSigma_{S}\|_{F}\leq 2\|\vD\|_{2}\|J_{S}\vV-\vV\|_{F}.
\end{align*}
Furthermore, if $\|J_{S}\vV-\vV\|_{F} \leq \frac{1}{8\kappa}$, then  $\|\bSigma-\bSigma_{S}\|_{F}\leq \lambda_{d}/4$. By Sin-Theta Lemma ( e.g., Lemma \ref{lem:sin_theta} ), we have
\begin{align*}
\|\vV\vV^{\tau}-\vV_{1}\vV_{1}^{\tau}\|_{F}\leq 8\kappa\|J_{S}\vV-\vV\|_{F}.
\end{align*}
\proof This comes from a (trivial) elementary calculus.  \epf
%\begin{align}
%\|\vV\vV^{\tau}-\vV_{1}\vV_{1}^{\tau}\|^{2}_{F}\leq \frac{2\|\bSigma-\bSigma_{S}\|^{2}_{F}}{\lambda^{2}_{m}(\Sigma_{S})}\leq 4\|\vD\|_{2}\|\vD^{-1}_{1}\|^{-1}_{2}\|J_{S}\vV-\vV\|^{2}_{F}.
%\end{align}
\end{lem}

Second, we need to introduce some notations and the `Oracle estimate' $\widehat{\vV}_{O}$. Since we have randomly divided the samples into two equal sets of samples, we have the corresponding decomposition  \eqref{eqn:estimator:decomposition}
\begin{equation*}
\begin{aligned} 
\widehat{\bLambda}_{H}
&=\bLambda_{u}+\bold{\mathcal{Z}}\bold{\mathcal{E}}^{\tau}\vV^{\perp,\tau}
+\vV^{\perp}\bold{\mathcal{E}}\bold{\mathcal{Z}}^{\tau}
+\vV^{\perp}\bold{\mathcal{E}}\bold{\mathcal{E}}^{\tau}\vV^{\perp,\tau}.
\end{aligned}
\end{equation*}
for these two sets of samples.  
More precisely, 
for $i=1,2,$ we can define $\bLambda_{H}^{(i)}$, $\bLambda_{\vu}^{(i)}$,  $\mathcal{Z}^{(i)}$, $\mathcal{W}^{(i)}$ and $\mathcal{E}^{(i)}$  for the first and second set of samples respectively according to the decomposition \eqref{eqn:estimator:decomposition}.
Let 
$\bLambda=\vV\vD\vV^{\tau}$ be the spectral decomposition,  
where  
$\vV$ is $p\times d$ orthogonal matrix and $\vD=diag\{\lambda_{1},..,\lambda_{d}\}$ is a diagonal matrix.
For $i=1,2,$  let
$\bLambda^{(i)}_{\vu}=\vV^{(i)}\vD^{(i)}\vV^{(i),\tau}$ ,
where  
$\vV^{(i)}$ is $p\times d$ orthogonal matrix and $\vD^{(i)}=diag\{\lambda_{1}^{(i)},..,\lambda_{d}^{(i)}\}$ is a diagonal matrix.
For any subset $S$ of $[p]$, let $J_{S}$ be the diagonal matrix  defined in Lemma \ref{sparse:random}.
Let $J_{S}\bLambda J_{S}=\vV_{S}\vD_{S}\vV^{\tau}_{S}$ be the spectral decomposition, where  
$\vV_{S}$ is $p\times d$ orthogonal matrix and $\vD_{S}=diag\{\lambda_{1,S},..,\lambda_{d,S}\}$ is a diagonal matrix.
Let $J_{S}\bLambda^{(i)}_{u}J_{S}= \vV^{(i)}_{S}\vD_{S}^{(i)}\vV^{(i),\tau}_{S}$ be the spectral decomposition, where  
$\vV_{S}^{(i)}$ is $p\times d$ orthogonal matrix and $\vD_{S}^{(i)}=diag\{\lambda_{1,S}^{(i)},..,\lambda_{d,S}^{(i)}\}$ is a diagonal matrix.
In the below, we will call $\vV_{S}$ ( {\it resp.}  $\vV_{S}^{(i)}$, $i=1,2$) the sparse approximation of  $\vV$ ( {\it resp.}  $\vV^{(i)}$, $i=1,2$).
From now on, we will choose $S$ to be $[k_{q,s}]\subset[p]$ where $k_{q,s}$ is defined near \eqref{def:sparse orthogonal}.

\vspace*{3mm}
%In the below, an event $\Omega$ is called happens with high probability if $\bbP(\Omega)\geq $
Below, we use $C$ to denote an absolute constant, though its exact value may vary from case to case.
We also assume that $\epsilon_{n}^{2}$ is sufficiently small. 
For $i=1,2$, let $\ttE^{(i)}_{3}$ be the event defined similarly as  $\ttE_{3}$ (which is introduced near Corollary \ref{cor:minimax:key}). Conditioning on $\ttE=\ttE^{(1)}_{3}\cap \ttE^{(2)}_{3}$,
Lemma \ref{lem:minimx:key}  implies that
\begin{align}\label{fact:key}
(1-\frac{\kappa}{2\nu})\lambda \leq \lambda^{(i)}_{d}\leq ... \leq \lambda^{(i)}_{1}\leq (1+\frac{1}{2\nu})\kappa\lambda \mbox{\quad \quad for } i=1,2.
\end{align}
%If we are working with the exact sparsity, then we have $\widetilde{\vV}^{(i)}_{S}=\widetilde{\vV}^{(i)}$. 
%If we are working with the weak-$l_{q}$ sparsity, we have the following results regarding the sparse approximation.
We first prove the following sparse approximation lemma.
\begin{lem}\label{true:sparse:approximation} Conditioning on $\ttE$, we have
\begin{align}
%\|%\widetilde{\vV}_{S}^{(i)}\vD_{H}^{(2)}\widetilde{\vV}_{S}^{(i),%\tau}-\widetilde{\vV}^{(i)}\vD_{H}^{(2)}\widetilde{\vV}^{(i),%\tau}\|_{F} &\leq C\sqrt{\frac{q}{2-%q}}\lambda\epsilon_{n} %\label{approximation:sparse:random}\\
\|\vV_{S}^{(i)}\vV_{S}^{(i),\tau}-\vV_{S}\vV_{S}^{\tau}\|^{2}_{F} &\leq C\frac{q}{2-q}\epsilon_{n}^{2} 
, \quad 
\|J_{S}\bLambda^{(i)}_{u}J_{S}-\bLambda^{(i)}_{u}\|_{F}
\leq C\lambda\epsilon_{n}\label{sparsebound:random}
\end{align}

and the entries of $\vD_{S}^{(i)} \in (\frac{1}{2}\lambda,2\kappa\lambda)$ for $i=1,2$.
\proof 
Since $\vV^{(i)}$ and $\vV$ share the same column space, we have $\vV^{(i)}=\vV\tilde{U}$ for some (stochastic) orthogonal matrix $\tilde{U}$ and $\vV\vV^{\tau}=\vV^{(i)}\vV^{(i),\tau}$. 
From this we know that %$\vV^{(i)}$ and $\vV$ are sharing the same sparse approximation set $S$ and 
\begin{align}
\nonumber &\|\vV_{S}^{(i)}\vV_{S}^{(i),\tau}-\vV_{S}\vV_{S}^{\tau}\|_{F} \leq 
\|\vV\vV^{\tau}-\vV_{S}\vV_{S}^{\tau}\|_{F} 
+\|\vV_{S}^{(i)}\vV_{S}^{(i),\tau}-\vV^{(i)}\vV^{(i),\tau}\|_{F}.
\end{align}

Conditioning on $\ttE$,  Lemma \ref{lem:risk:sparse:approximation}, Lemma \ref{sparse:random} and \eqref{fact:key} imply 
%\begin{equation}
%\begin{aligned}
%\nonumber
%\|J_{S}\bLambda^{(i)}_{u}J_{S}-\bLambda^{(i)}_{u}\|_{F}
%\leq 2(1+\frac{1}{2\nu})\kappa\lambda\|J_{S}\widetilde{\vV}^{(i)}-\widetilde{\vV}^{(i)}\|
%\leq 2(1+\frac{1}{2\nu})\kappa\lambda\epsilon_{n}.%\label{sparsebound:random}
%\end{aligned}
%\end{equation}
\begin{equation*}
\|J_{S}\bLambda^{(i)}_{u}J_{S}-\bLambda^{(i)}_{u}\|_{F}
\leq C\kappa\lambda\|J_{S}\vV^{(i)}-\vV^{(i)}\|_{F}
\leq C\kappa\lambda\epsilon_{n}.
\end{equation*}
Since we have assumed that $\epsilon^{2}_{n}$ is sufficiently small,
we can assure that the entries of $\vD_{S}^{(i)}$ are in the range $(\frac{1}{2}\lambda,2\kappa\lambda)$.
After applying the Sin-Theta theorem (e.g. Lemma \ref{lem:sin_theta}), we have
\begin{align*}
\|\vV_{S}^{(i)}\vV_{S}^{(i),\tau}-\vV^{(i)}\vV^{(i),\tau}\|_{F}\leq C\kappa\|J_{S}\vV^{(i)}-\vV^{(i)}\|_{F}\leq C\kappa\sqrt{\frac{q}{2-q}}\epsilon_{n}.
\end{align*}
We can apply similar argument to bound $\|\vV\vV^{\tau}-\vV_{S}\vV_{S}^{\tau}\|_{F}$, which gives us
\begin{align}
\nonumber &\|\vV_{S}^{(i)}\vV_{S}^{(i),\tau}-\vV_{S}\vV_{S}^{\tau}\|_{F}\leq C\kappa\|J_{S}\vV-\vV\|_{F}\leq C\sqrt{\frac{q}{2-q}}\epsilon_{n}.
\end{align} \epf
\end{lem}

We introduce an `Oracle estimator' $\widehat{\vV}_{O}$ (as if we know the sparse approximation set $S$)  such that
\begin{equation}\label{estimator:oracle}
\begin{aligned}
&\widehat{\vV}_{O}\triangleq \arg\max_{\vV} \langle\bLambda^{(1)}_{H},\vV\vV^{\tau}\rangle
=\arg\max_{\vV}Tr(\vV^{\tau}\bLambda^{(1)}_{H}\vV)\\
&\mbox{ s.t. }  \vV^{\tau}\vV=I_{d} \mbox{ and } supp(\vV) = S.
\end{aligned}
\end{equation}  
Let $\widehat{\vV}_{O}^{\tau}\vV_{S}^{(2)}=U_{1}\Delta U_{2}^{\tau}$ be the singular value decomposition of $\widehat{\vV}_{O}^{\tau}\vV_{S}^{(2)}$ such that the entries of $\Delta$ are non-negative and let $M\triangleq U_{2}^{\tau}\vD_{S}^{(2)}U_{2}$.  

\vspace*{6mm} Now, we can start our proof of Theorem \ref{thm:risk:sparse:upper:d}. It is easy to verify that
\begin{align*}
\|\widehat{\vV}_{E}\widehat{\vV}_{E}^{\tau}-\vV\vV^{\tau}\|_{F}^{2}
\leq C\left( 
\|\widehat{\vV}_{E}\widehat{\vV}_{E}^{\tau}-\widehat{\vV}_{O}\widehat{\vV}_{O}^{\tau}\|_{F}^{2}
+\|\widehat{\vV}_{O}\widehat{\vV}_{O}^{\tau}-\vV_{S}\vV_{S}^{\tau}\|_{F}^{2}+
\|\vV_{S}\vV_{S}^{\tau}-\vV\vV^{\tau}\|_{F}^{2} \right).
\end{align*}
For the first term $\|\widehat{\vV}_{E}\widehat{\vV}_{E}^{\tau}-\widehat{\vV}_{O}\widehat{\vV}_{O}^{\tau}\|^{2}_{F}$, conditioning on $\ttE$, we know 
\begin{align}
\|\widehat{\vV}_{E}\widehat{\vV}_{E}^{\tau}-\widehat{\vV}_{O}\widehat{\vV}_{O}^{\tau}\|^{2}_{F} 
\leq & \frac{2}{\lambda_{d}(\vD_{S}^{(2)})}\langle\widehat{\vV}_{O}U_{1}MU_{1}^{\tau}\widehat{\vV}_{O}^{\tau},\widehat{\vV}_{O}\widehat{\vV}_{O}^{\tau}-\widehat{\vV}_{E}\widehat{\vV}_{E}^{\tau} \rangle \label{inline:E-1:temp}\\
\leq & \frac{C}{\lambda}\langle\widehat{\vV}_{O}U_{1}MU_{1}^{\tau}\widehat{\vV}_{O}^{\tau}-\bLambda^{(2)}_{H},\widehat{\vV}_{O}\widehat{\vV}_{O}^{\tau}-\widehat{\vV}_{E}\widehat{\vV}_{E}^{\tau} \rangle \label{inline:E-O:temp} \\
\nonumber \triangleq& I+II+III.
\end{align}
where
\begin{align}
I\nonumber=&\frac{C}{\lambda}\langle\widehat{\vV}_{O}U_{1}MU_{1}^{\tau}\widehat{\vV}_{O}^{\tau}
-\vV^{(2)}_{S}\vD_{S}^{(2)}\vV_{S}^{(2),\tau},\widehat{\vV}_{O}\widehat{\vV}_{O}^{\tau}-\widehat{\vV}_{E}\widehat{\vV}_{E}^{\tau} \rangle \\
II\nonumber=&\frac{C}{\lambda}\langle \vV^{(2)}_{S}\vD_{S}^{(2)}\vV_{S}^{(2),\tau}-\bLambda^{(2)}_{\vu},\widehat{\vV}_{O}\widehat{\vV}_{O}^{\tau}-\widehat{\vV}_{E}\widehat{\vV}_{E}^{\tau} \rangle & \\
III\nonumber =&\frac{C}{\lambda}\langle\bLambda^{(2)}_{\vu}-\bLambda^{(2)}_{H},\widehat{\vV}_{O}\widehat{\vV}_{O}^{\tau}-\widehat{\vV}_{E}\widehat{\vV}_{E}^{\tau} \rangle. &
\end{align}
Inequality \eqref{inline:E-1:temp} follows from applying the Lemma \ref{cor:elemetary:temp} with the positive definite matrix $U_{1}MU_{1}^{\tau}$.  
The inequality \eqref{inline:E-O:temp} follows from the definition of $\widehat{\vV}_{E}$ (See \eqref{eqn:two:set:inequality}) and the fact that the entries of $\vD_{S}^{(2)}$ are in $(\lambda/2,2\kappa\lambda)$. To simplify the notation, we let
\begin{equation*}
\begin{aligned}
R=\|\widehat{\vV}_{E}\widehat{\vV}_{E}^{\tau}-& \vV\vV^{\tau}\|_{F},  \theta^{(i)}=\|\vV_{S}^{(i)}\vV_{S}^{(i),\tau}-\vV\vV^{\tau}\|_{F}, 
\delta=\|\widehat{\vV}_{O}\widehat{\vV}_{O}^{\tau}
-\vV_{S}^{(1)}\vV_{S}^{(1),\tau}\|_{F}.
\end{aligned}
\end{equation*}
%$\delta_{1}=\|\widehat{\vV}_{O}\widehat{\vV}_{O}^{\tau}-\vV_{S}\vV_{S}^{\tau}\|_{F}$
%and
% $\delta_{2}=\|\vV_{S}\vV_{S}^{\tau}-\widetilde{\vV}_{S}^{(2)}\widetilde{\vV}_{S}^{(2),\tau}\|_{F}$.

\paragraph{For I:} 
Recall that the entries of $\vD_{S}^{(2)}$ $\in (\frac{1}{2}\lambda, 2\kappa\lambda)$  and that $\widehat{\vV}_{O}U_{1}$ and $\vV_{S}^{(2)}U_{2}$   satisfies the condition that $U_{1}^{\tau}\widehat{\vV}_{O}^{\tau}\vV_{S}^{(2)}U_{2}$ is a diagonal matrix with non-negative entries.
Since $M\triangleq U_{2}^{\tau}\vD_{S}^{(2)}U_{2}$,  by Lemma \ref{lem:elementary:trivial3} and the fact that $\vV_{S}^{(2)}U_{2}$  and  $\vV_{S}^{(2)}$ share the same column space,  there exists a constant $C$ such that
%\begin{equation}
%\begin{aligned}
%\|\widehat{\vV}_{O}U_{1}MU_{1}^{\tau}\widehat{\vV}_{O}^{\tau}-\bLambda_{\vu}^{(2)}\|_{F}&= \|\widehat{\vV}_{O}U_{1}MU_{1}^{\tau}\widehat{\vV}_{O}^{\tau}-\widetilde{\vV}^{(2)}
%U_{2}MU_{2}^{\tau}\widetilde{\vV}^{(2),\tau}\|_{F}\\
%&\leq C\lambda \|\widehat{\vV}_{O}\vV^{\tau}_{O}-\widetilde{\vV}^{(2)}\widetilde{\vV}^{(2),\tau}\|_{F}\\
%&=C\lambda \|\widehat{\vV}_{O}\vV^{\tau}_{O}-\vV\vV^{\tau}\|_{F},
%\end{aligned}
%\end{equation}
\begin{equation}\nonumber
\begin{aligned}
\| \widehat{\vV}_{O}U_{1}MU_{1}^{\tau}\widehat{\vV}_{O}^{\tau}
-\vV^{(2)}_{S}\vD_{S}^{(2)}\vV^{(2),\tau}_{S} \|_{F}
\leq C\lambda \|\widehat{\vV}_{O}\widehat{\vV}^{\tau}_{O}-\vV_{S}^{(2)}\vV_{S}^{(2),\tau}\|_{F}.
\end{aligned}
\end{equation}
Thus, conditioning on $\ttE$, we have
\begin{equation}\label{inlince:ttt}
\begin{aligned}
\big|I\big| \leq& C \|\widehat{\vV}_{O}\widehat{\vV}_{O}^{\tau}-
\vV_{S}^{(2)}\vV_{S}^{(2),\tau}\|_{F}\|\widehat{\vV}_{O}\widehat{\vV}_{O}^{\tau}-\widehat{\vV}_{E}\widehat{\vV}_{E}^{\tau}\|_{F}\\
\leq & C ( \delta+\theta^{(1)}+\theta^{(2)} )\|\widehat{\vV}_{O}\widehat{\vV}_{O}^{\tau}-\widehat{\vV}_{E}\widehat{\vV}_{E}^{\tau}\|_{F}.
\end{aligned}
\end{equation}

\paragraph{For II:}  It is nonzero only if $q\neq 0$.
From  Lemma \ref{lem:risk:sparse:approximation} and Lemma \ref{sparse:random}, we know that 
\begin{align}
|II|\leq C\sqrt{\frac{q}{2-q}}\epsilon_{n}\|\widehat{\vV}_{O}\widehat{\vV}_{O}^{\tau}-\widehat{\vV}_{E}\widehat{\vV}_{E}^{\tau}\|_{F}.
\end{align}

\paragraph{For III:}  
From the equation \eqref{eqn:estimator:decomposition}, we have
\begin{equation}\label{inlince:tt}
\begin{aligned}
\big|III \big| 
\leq& \frac{1}{\lambda}\|\widehat{\vV}_{O}\widehat{\vV}_{O}^{\tau}-\widehat{\vV}_{E}\widehat{\vV}_{E}^{\tau}\|_{F}\left( 2T_{2}+T_{1}\right)
\end{aligned}
\end{equation}
where 
$T_{1}=\max_{B\in\mathcal{B}(k)}\Big|\Big< \mathcal{W}^{(2)}\mathcal{W}^{(2),\tau},\vK_{B}\Big>\Big|$,
$T_{2}=\max_{B\in\mathcal{B}(k)}\Big|\Big< \mathcal{Z}^{(2)}\mathcal{W}^{(2),\tau},\vK_{B}\Big>\Big|$ and $\vK_{B}=\|\widehat{\vV}_{O}\widehat{\vV}_{O}^{\tau}-\widehat{\vV}_{B}\widehat{\vV}_{B}^{\tau}\|^{-1}_{F}\left(\widehat{\vV}_{O}\widehat{\vV}_{O}^{\tau}-\widehat{\vV}_{B}\widehat{\vV}_{B}^{\tau}\right)$. ( For any $B \in \mathcal{B}_{k}$, $\widehat{\vV}_{B}$ is introduced in \eqref{estimator:sample} ).

\vspace*{3mm}
To summarize, conditioning on E, we have
\begin{align}
\|\widehat{\vV}_{E}\widehat{\vV}_{E}^{\tau}-
\widehat{\vV}_{O}\widehat{\vV}_{O}^{\tau}\|_{F}\leq C\left(\delta+\theta^{(1)}+\theta^{(2)}+\epsilon_{n}+\frac{1}{\lambda}(2T_{2}+T_{1})\right).
\end{align}
Thus, we have
\begin{align*}
R^{2}\bold{1}_{\ttE} 
&\leq 
C\left( \delta^2+\left(\theta^{(1)}\right)^{2}+\|\widehat{\vV}_{O}\widehat{\vV}_{O}^{\tau}-\widehat{\vV}_{E}\vV^{\tau}_{E}\|^{2}_{F}\right)\bold{1}_{\ttE}\\
&\leq 
C\left( \delta^2+\left(\theta^{(1)}\right)^{2}+C\left( \delta+\theta^{(1)}+\theta^{(2)} +\epsilon_{n}+\frac{1}{\lambda}\left(2T_{1}+T_{2} \right)\right)^2\right)\bold{1}_{\ttE}\\
&\leq C\left( \delta^2+\left(\theta^{(1)}\right)^{2}+\left(\theta^{(2)}\right)^{2} +\epsilon^{2}_{n}\bold{1}_{\ttE}+\left(\frac{1}{\lambda}\left(2T_{1}+T_{2} \right)\right)^2\right)\bold{1}_{\ttE}
\end{align*}
If we can prove
\begin{align}
\bbE\left( \theta^{(i)}\right)^{2}\bold{1}_{\ttE'}\leq C\epsilon_{n}^{2}, \quad
\bbE\delta^{2}\bold{1}_{\ttE'}\leq C\epsilon_{n}^{2} \quad  
\mbox{ and } \bbE( 2T_{1}+T_{2})^{2}\bold{1}_{\ttE'}\leq \lambda^{2} \epsilon^{2}_{n},
\end{align}   
for some $\ttE'\subset \ttE$ such that$\bbP\left(\left(\ttE'\right)^{c}\right)\leq C\frac{H^{2}}{n\lambda}$,
then we have $\bbE R^{2}\bold{1}_{\ttE'}\leq C\epsilon_{n}^2$. 
%Note that $\bbP(\ttE^{c})\leq C\frac{H^{2}}{n\lambda}$ (See Corollary \ref{cor:minimax:key}) implies $\bbE R^{2}\bold{1}_{\ttE^{c}} \leq C \frac{dH^{2}}{n\lambda}$. 
Thus, we have
\[
\bbE R^{2} \leq C\epsilon_{n}^{2}.
\]\epf

All we need to  prove are the following two Lemmas.
\begin{lem} Assume that $n\lambda\leq e^{p}$. There exist $\ttE'\subset \ttE$ such that  $\bbP((\ttE')^{c})\leq C\frac{H^{2}}{n\lambda}$ and
\begin{align}
\bbE\left( \theta^{(i)}\right)^{2}\bold{1}_{\ttE'}\leq C\epsilon_{n}^{2} \mbox{ \quad and \quad }
\bbE\delta^{2}\bold{1}_{\ttE'}\leq C\epsilon_{n}^{2}.  
\end{align}
\proof 
Since $\vV^{(i)}$ and $\vV$ share the same column space, conditioning on $\ttE$, by Lemma \ref{sparse:random} and  Lemma \ref{lem:risk:sparse:approximation}, we have
\begin{equation}  \label{partone:dec}
\begin{aligned}  
\theta^{(i)}=\|\vV_{S}^{(i)}\vV_{S}^{(i),\tau}-\vV\vV^{\tau}\|_{F}\leq 4\kappa\|J_{S}\vV^{(i)}-\vV^{(i)}\|_{F}\leq C\sqrt{\frac{q}{2-q}}\epsilon_{n},
\end{aligned}
\end{equation}
i.e. $\bbE\left( \theta^{(i)}\right)^{2}\bold{1}_{\ttE}\leq C\epsilon_{n}^{2}$.

Let $\vQ_{S}= J_{S} \left(\bLambda^{(1)}_{H}-\bLambda^{(1)}_{u}\right)J_{S}$.
Let $\mathtt{F}$ consist of the events such that $\|J_{S}\mathcal{\vW}^{(1)}\mathcal{\vW}^{(1),\tau}J_{S}\|\leq C\frac{k}{n}$. Lemma \ref{random:nonasymptotic} implies that $\bbP(\mathtt{F}^{c})\leq C\frac{H^{2}}{n\lambda}$.
%$\geq {\color{red} 1-2\exp ( -k\log \frac{ep}{k}) }$.
Since we have assumed that $\epsilon_{n}^{2}$ is sufficiently small, conditioning on $\ttE\cap\mathtt{F}$, the decomposition 
\eqref{eqn:estimator:decomposition}
%to the matrix $\sqrt{n}J_{S_{q,s}}\mathcal{W}$ with $t^2=\log \frac{ep}{k}$, we have
give us
%\begin{equation}\label{bound:sparse:noisy}
%\begin{aligned}
%\|J_{S_{q,s}} \bLambda^{(1)}_{H}J_{S_{q,s}}-J_{S_{q,s}}\bLambda^{(1)}_{u}J_{S_{q,s}}\| \leq \frac{1}{n}\log \frac{ep}{k}+2\sqrt{\frac{1}{n}\log \frac{ep}{k}\lambda (1+\frac{1}{2\nu})}\leq C\lambda \sqrt{ \frac{1}{n\lambda}\log \frac{ep}{k} },
%\end{aligned}
%\end{equation}
\begin{align}\label{bound:sparse:noisy}
\|\vQ_{S}\|_{2}\leq C\sqrt{\lambda}\sqrt{\frac{k}{n}}\leq C\frac{\lambda\epsilon_{n}}{\sqrt{d}}.
\end{align} 
%From Corollary \ref{cor:minimax:key}, we have 
%\begin{equation}\label{bound:sparse:noisy}
%\begin{aligned}
%\|J_{S_{q,s}}\bLambda^{(1)}_{H}J_{S_{q,s}}-J_{S_{q,s}}\Lambda^{(1)}_{u}J_{S_{q,s}}\|_{}\leq C\lambda \sqrt{\kappa  \frac{k_{q,s}}{n\lambda}}.
%\end{aligned}
%\end{equation} 
%Assuming $\epsilon^{2}_{n}$ is sufficiently small, the bound above is smaller than $\frac{\lambda}{8}$. 
From \eqref{sparsebound:random}, we also have $\|J_{S}\bLambda^{(1)}_{u}J_{S}-\bLambda^{(1)}_{u}\|_{F}\leq C\lambda\epsilon_{n}\leq \frac{\lambda}{8}$. Thus, $\|J_{S}\bLambda^{(1)}_{H}J_{S} - \bLambda^{(1)}_{u}\|_{F}< \frac{\lambda}{4} $, which implies the $(d+1)$-th largest eigenvalues of $J_{S}\bLambda^{(1)}_{H}J_{S}$ is less than $\frac{\lambda}{4}$. 
Note that the eigenvalues of $J_{S}\bLambda^{(1)}_{u}J_{S}\in (\frac{1}{2}\lambda,2\kappa\lambda)$. After applying the Sin-Theta Theorem( Lemma \ref{lem:sin_theta}) to  the pair of symmetric matrices $(J_{S}\bLambda^{(1)}_{u}J_{S}, J_{S}\bLambda^{(1)}_{H}J_{S})$ , we have
\begin{equation*}\label{use:sintheta}
\begin{aligned}
\delta 
\leq  \frac{8}{\lambda} \| \widehat{\vV}_{O}^{\perp,\tau } \vQ_{S}  \vV_{S}^{(1)}\|_{F}
\leq  \frac{8}{\lambda} \sqrt{d}\|\vQ_{S}\|_{2} 
\leq C\epsilon_n
 \end{aligned}
\end{equation*} 
where the last inequality follows from \eqref{bound:sparse:noisy}. Thus, we may take $\ttE'=\ttE\cap\mathtt{F}$.
\end{lem}

\begin{lem} There exists positive constant $C$ such that
\[
\bbE(2T_{1}+T_{2})^{2} \bold{1}_{\ttE}\leq C\lambda^{2}\epsilon^{2}_{n}
\]
\proof Since $(2T_{1}+T_{2})^{2} \leq C(T_{1}^{2}+T_{2}^{2})$, we only need to bound $\bbE T_{1}^{2}$ and $\bbE T_{2}^{2}$ separately.
\paragraph{For $T_{1}$}
Recall that $\mathcal{W}^{(2)}=\vV^{\perp}\mathcal{E}^{(2)}$ (See notation near \eqref{eqn:estimator:decomposition}.) and for each fixed $B \in \mathcal{B}_{k}$, $\vK_{B}\independent \mathcal{W}^{(2)}$, hence 
\begin{align}
\langle\mathcal{W}^{(2)}\mathcal{W}^{(2),\tau},\vK_{B} \rangle=\langle\mathcal{E}^{(2)}\mathcal{E}^{(2),\tau},\vV^{\perp,\tau}\vK_{B}\vV^{\perp} \rangle
\end{align} 
and $\vV^{\perp,\tau}\vK_{B}\vV^{\perp}  \independent \mathcal{W}^{(2)}$. 
Note that $\|\vV^{\perp,\tau}\vK_{B}\vV^{\perp}\|_{F}\leq 1$, $\mathcal{E}^{(2)}$ is a $(p-d)\times H$ matrix and $\sqrt{n}\mathcal{E}_{i,j}^{(2)}\sim N(0,1)$. After applying Lemma \ref{lem:Cai:lem4},  we have
\begin{align}
\bbP\left(\sqrt{n}|\langle \mathcal{E}^{(2)}\mathcal{E}^{(2),\tau},\vV^{\perp,\tau}\vK_{B}\vV^{\perp} \rangle| \geq 2\frac{\sqrt{H}}{\sqrt{n}} t+\frac{2}{\sqrt{n}}t^{2} \right) \leq 2\exp	(-t^{2}).
\end{align}
After applying Lemma \ref{lem:Cai:lem5} with $N=|\mathcal{B}(k)|\leq \left(\frac{ep}{k} \right)^{k}$, $a=\frac{2\sqrt{H}}{\sqrt{n}}$, $b=\frac{2}{\sqrt{n}}$ and $c=2$, we have
\begin{align*}
\bbE T_{1}^{2} & \leq \frac{1}{n} \left( \frac{8H}{n}\log(2eN) +\frac{8}{n}\left(\log^{2}(2N)+4\log(2eN) \right) \right)\\
&=\frac{8(H+4)\log(2eN)}{n^{2}}+\frac{8}{n^{2}}\log^{2}(2N)\\
&\leq C\lambda^{2}\frac{\log(N)}{n\lambda}\leq C\lambda^{2}\epsilon_{n}^{2}
\end{align*}

\paragraph{For $T_{2}$}
 Fix $B \in \mathcal{B}(k_{q,s})$. Since $\mathcal{Z}^{(2)}\independent \mathcal{W}^{(2)}$, $\vK_{B}\independent \mathcal{W}^{(2)}$ and  $\vK_{B}\independent \mathcal{Z}^{(2)}$, conditioned on the $\mathcal{Z}^{(2)}$ and $\vK_{B}$,  we know that 
\[
\sqrt{n}\langle \mathcal{Z}^{(2)}\mathcal{W}^{(2),\tau},\vK_{B}\rangle=\langle\vV^{\perp,\tau}\vK_{B}\mathcal{Z}^{(2)},\sqrt{n}\mathcal{E}^{(2)} \rangle
\] is distributed according to $N(0, \|\vV^{\perp,\tau}\vK_{B}\mathcal{Z}\|^{2}_{F})$. 
Therefore 
\begin{align*}
\sqrt{n}\langle \mathcal{Z}^{(2)}\mathcal{W}^{(2),\tau},\vK_{B}\rangle\overset{d}{=}\|\vV^{\perp,\tau}\vK_{B}\mathcal{Z}^{(2)}\|_{F}W
\end{align*}
for some $W \sim N(0,1)$ independent of $\mathcal{Z}^{(2)}$ and $\vK_{B}$. For simplicity of notation, we denote $ \sqrt{n}\langle \mathcal{Z}^{(2)}\mathcal{W}^{(2),\tau},\vK_{B} \rangle$ by $F_{B}$. As a direct corollary,  conditioning on $\ttE$, we know
\begin{align}
\bbP\left(|F_{B}| >t\right) \leq \bbP\left(2\kappa\lambda |W| >t\right) \leq 2 \exp\left(-\frac{t^{2}}{4\kappa^{2}\lambda^{2}} \right).
\end{align}
i.e., conditioning on $\ttE$, $|F_{B}|$ is sub-Gaussian and upper exponentially bounded by $4\kappa^{2}\lambda^{2}$. From this, we know $\bbE \left( T_{2}^{2} \bold{1}_{\ttE}\right)\leq C\lambda^{2}\epsilon^{2}_{n}$.

\epf
\end{lem}

\iffalse
 % and the fact that Frobenius norm dominates operator norm.
\begin{equation}
\begin{aligned}
 &\|\widehat{\vV}_{O}\vV^{\tau}_{O}-\widetilde{\vV}_{S}^{(2)}\widetilde{\vV}_{S}^{(2),\tau}\|_{F}\\
 \leq& \|\widehat{\vV}_{O}\vV^{\tau}_{O}-\widetilde{\vV}_{S}^{(1)}\widetilde{\vV}_{S}^{(1),\tau}\|_{F}+\| \widetilde{\vV}_{S}^{(1)}\widetilde{\vV}_{S}^{(1),\tau} - \widetilde{\vV}^{(1)}\widetilde{\vV}^{(1),\tau} \|_{F}\\
 &+
\|\widetilde{\vV}^{(2)}\widetilde{\vV}^{(2),\tau} - \widetilde{\vV}_{S}^{(2)}\widetilde{\vV}_{S}^{(2),\tau}\|_{F}\\
 =&\delta+\theta^{(1)}+\theta^{(2)}
\end{aligned}
\end{equation}
\fi

\section{Discussion}\label{sec:conclusion}
In this paper, we have determined the minimax rate of estimating the central space over a large class  of models $\mathfrak{M}_{s,q}\left( p,d,\lambda,\kappa\right)$ in two scenarios: $1)$ single index models and $2)$ $d$ and $\lambda$ are bounded. Here $\lambda$, the smallest nonzero eigenvalue %$\lambda_{\min}(var(\bbE[\vx|y]))$ 
of $var(\bbE[\vx|y])$, plays the role of signal strength in SIR and can be viewed as a generalized notion of the signal-to-noise ratio for multiple index models. Since we have established an upper bound of convergence rate of estimating the central space for all $d$ and $\lambda$, we will attempt to show that this convergence rate is optimal even for diverging $d$ and $\lambda$ in a future research. 
%{\color{green}
%More precisely, if we consider the simplest linear model, $y=\bbeta^{\tau}\vx+\epsilon$, where $\vx\sim N(0,\vI)$ and $\epsilon\sim N(0,\sigma^{2})$, then $\lambda_{\min}(var(\bbE[\vx|y]))=\frac{\|\bbeta\|^{2}_{2}}{\|\bbeta\|^{2}_{2}+\sigma^{2}}$ is the signal noise ratio (SNR) in the linear model. In other words, we could treat the smallest nonzero eigenvalue of $var(\bbE[\vx|y])$ as a generalized notation of SNR for multiple index models. }

The aggregate estimator we constructed here is actually an estimator of the column space of $var(\bbE[\vx|y])$ rather than that of the central space. Since we have assumed that $\bSigma=\vI$ in this paper, the column space of $var(\bbE[\vx|y])$ coincides with the central space in model \eqref{model:multiple}. 
When there are correlations between predictors,
if we assume that the eigenvectors associated with non-zero eigenvalues of $var(\bbE[\vx|y])$ are sparse (with sparsity $s$) instead of assuming that the loading vectors $\bbeta_{i}'s$ are sparse, our argument in this paper implies  that $\bbE[\|P_{\widehat{col(var(\bbE[\vx|y]))}}-P_{col(var(\bbE[\vx|y]))}\|_{F}^{2}]$ converges at the rate $\frac{ds+s\log(ep/s)}{n\lambda}$.

Although our studies of the sparse SIR were inspired by  recent advances in sparse PCA,  the results in this paper suggest a more intimate connection between SIR and linear regressions. Recall that
%  a more appropriate prototype of SIR in high dimensions might be the sparse linear regressions, at least for single index models. 
for the linear regression model $y=\bbeta^{\tau}\vx+\epsilon$ with $\vx\sim N(0,\vI)$ and $s=O(p^{1-\delta})$, 
the minimax rate  \citep{raskutti2011minimax} of estimating $\bbeta$  is achieved by the simple correlation screening. On the other hand, 
the  minimax rate for estimating $P_{\bbeta}$ is achieved by the DT-SIR algorithm of \cite{lin2015consistency}, which simply screens each variable based on the estimated variance of its conditional means. This fact suggests that a more appropriate prototype of SIR in high dimensions might be linear regression rather than sparse PCA, because 
there is a computational barrier of the rate optimal estimates for sparse PCA \citep{berthet2013computational}.
This possibility further suggests  that an efficient (rate optimal) high dimensional variant of SIR with general variance matrix $\bSigma$ might be possible, since it is now well known that Lasso\citep{tibshirani1996regression} and Dantzig Selector\citep{candes2007dantzig} achieve the optimal rate of linear regression \citep{bickel2009simultaneous} for general $\bSigma$. This speculation warrants further future investigations.

\begin{supplement}\label{suppA}
%\sname{Supplement}
\stitle{Supplement to `` On the optimality of  SIR in high dimensions''}
\slink[url]{http://www.e-publications.org/ims/support/dowload/imsart-ims.zip}
\sdescription{}
\end{supplement}

\bibliographystyle{plainnat}
\bibliography{sir}

\newpage
\begin{frontmatter}

\title{Supplement to : ``On optimality of sliced inverse regression  in high dimensions''}
\runtitle{Supp: DT-SIR }

\begin{aug}
 \author{\fnms{Qian} \snm{Lin}\thanksref{m1}\thanksref{}}
 \author{\fnms{Xinran} \snm{Li}\thanksref{m1}\thanksref{}}
  \author{\fnms{Dongming} \snm{Huang}\thanksref{m1}\thanksref{}}
\and    
\author{\fnms{Jun S.} \snm{Liu}\thanksref{m1}\thanksref{}}

\runauthor{Q. Lin, X. Li, D. Huang and J. S. Liu}

\affiliation{Harvard University\thanksmark{m1} }

%\affiliation{Temple University\thanksmark{m2} }
  
\end{aug}

\end{frontmatter}

\begin{appendix}

%\section{Additional proofs}

\subsection*{Proof of Theorem \ref{thm:DTSIR}}
For a vector $\gamma \in \bbR^{p}$, $S\subset [p]$, let $\gamma_{S} \in \bbR^{p}$ such that $\gamma_{S}(i)=\gamma(i)$ if $i \in S$ and  $\gamma_{S}(i)=\gamma(i)$ if $i \not \in S$. For any non-zero vector $\gamma$,	let $\widetilde{\gamma}=\gamma/\|\gamma\|_{2}$.
For any non-zero vector $\gamma$ and $t>0$, let $T_{t}$ be the indices such that $|\gamma(i)|>t$.
We have following elementary Lemmas. 
\begin{lem}\label{lem:temp:supp} Let $\bgamma\in \mathbb{R}^{p}$ be a unit vector with at most $s$ non-zero entries, then
\begin{align}
\|\bgamma-\widetilde{\bgamma}_{T_{t}}\|^{2}_{2}\leq Cst^{2}.
\end{align}
\proof Let $E^{2}=\sum_{|\gamma_{j}^{}|\geq t} \gamma_{j}^{2}$, then
$
\|\gamma-\tau_{N}(\gamma,t)\|_{2}^{2}=\sum_{\gamma_{j}\geq t} \gamma_{j}^{2}(1-\frac{1}{E})^{2}+\sum_{\gamma_{j}<t}\gamma_{j}^{2}=2(1-E)\leq 2st^{2}.
$
\epf
\end{lem}

\begin{lem}
Let $\bLambda=\lambda\bbeta\bbeta^{\tau}$. For $S\subset [p]$, we have $\lambda\left(\bLambda(S,S)\right)=\lambda\|\bbeta_{S}\|^{2}\widetilde{\bbeta}_{S}\widetilde{\bbeta}_{S}^{\tau}$.
\proof It follows from (trivial) elementary calculus.\epf
\end{lem}

Let $T=\lbrace~i\mid \bLambda_{H}(i,i)>a\frac{\log(p)}{n} \rbrace$, then $\lambda_{T}=\lambda\left(\bLambda(T,T) \right)\geq \lambda\|\bbeta_{T}\|^{2}\geq \lambda(1-\frac{as\log(p)}{n})\geq a'\lambda$ if $\frac{s\log(p)}{n\lambda}$ is sufficiently small.
Since  $\widehat{\bLambda}_{H}(i,i)\sim \frac{1}{n}\chi^{2}_{H}$ for $i\not \in S$, we have $\bbP\left(\widehat{\bLambda}_{H}(i,i)>a\frac{\log(p)}{n}\right)\leq \exp\left(-a\log(p) \right)$. Thus, we have
$\bbP\left(T\subset S \right)\geq \bbP\left(\max_{i\not \in S} \widehat{\bLambda}_{H}(i,i)\leq a\frac{\log(p)}{n}\right)\geq 1-\exp(-(a-1)log(p))  $.

Thus, if $T\subset S$, we have 
\begin{align*}
\|\widehat{\bbeta}_{T}-\bbeta\|_{F}^{2}\leq \|\widehat{\bbeta}_{T}-\widetilde{\bbeta}_{T}\|_{F}^{2}+\|\widetilde{\bbeta}_{T}-\bbeta\|_{F}^{2}\leq C\frac{|T|}{n\lambda_{T}}+Cst^{2}\leq C\frac{s\log(p)}{n\lambda}
\end{align*}
where we have used the Oracle risk Theorem \ref{thm:risk:oracle:upper:d} and the Lemma \ref{lem:temp:supp}.
Thus, we know that DT-SIR is rate optimal if $s=O(p^{1-\delta})$.

\subsection*{The lower bound}
In this subsection, we provide the proof of the lower bound for 
%Theorem \ref{thm:risk:oracle:upper:d},
%Theorem\ref{thm:risk:sparse:upper:d},
Theorem \ref{thm:oracle:risk:d:fxied:lambda:fixed},
Theorem \ref{thm:rsik:sparse:d:fixed:lambda:fixed},
Theorem \ref{thm:risk:oracle:d=1}.
%Theorem\ref{thm:rsik:sparse:d=1},
%Theorem \ref{thm:oracle:risk},
%Theorem \ref{thm:sparse:rsik}.
%Theorem \ref{thm:DTSIR},
\paragraph{Proof of Theorem \ref{thm:oracle:risk:d:fxied:lambda:fixed}}
Let us consider the Grassmannian $\mathbb{G}(p,d)$ consisting of all the $d$ dimensional subspaces in $\bbR^{p}$ and the homogeneous space $\mathbb{O}(p,d)$ consisting of all $p\times d$ orthogonal matrices. There is a tautological map from $\mathbb{O}(p,d)$ to $\mathbb{G}(p,d)$, {\it i.e.}, $\vA\mapsto \vA\vA^{\tau}$.
For any $\varepsilon \in (0,\sqrt{2d\wedge(p-d)}]$, for any $u \in \mathbb{G}(p,d)$, \cite{cai2013sparse} have constructed a subset $\Theta \subset N(u,2\varepsilon)$, an $2\varepsilon$ neighbourhood of $u$ in $\mathbb{G}(p,d)$, such that, for any $\alpha \in (0,1)$, we have 
\begin{align*}
|u_{i}-u_{j}|\leq \varepsilon , \quad
|u_{i}-u_{j}| \geq \alpha \varepsilon \mbox{ and }
\left|\Theta\right|\geq \left(\frac{c_{0}}{\alpha c_{1}}\right)^{d(p-d)}
\end{align*}
where $u_{i}$ and $u_{j}$ are two different points $\in \Theta$ and $c_{0}$ and $c_{1}$ are two absolute constants.
Lemma \ref{lem:temp} states that if $\varepsilon$ is sufficiently small, then for each $u_{i} \in \Theta \subset \mathbb{G}(p,d)$, there is an $u_{i} \in \mathbb{O}(p,d)$ such that $a_{i}a_{i}^{\tau}=u_{i}$  
\begin{align}
C_{1}\|u_{i}-u_{j}\|_{F} \leq \|a_{i}-a_{j}\|_{F}\leq  C_{2}\|u_{i}-u_{j}\|_{F}
\end{align}
for some absolute positive constants $C_{1}$ and $C_{2}$. 
Let us denote $\widetilde{\Theta}=\{a_{i}\}$ and consider the following models
\begin{align}
y=f(\vV^{\tau}\vx)+\epsilon, \vV\in \widetilde{\Theta}, \vx \sim N(0,\vI_{p}), \epsilon\sim N(0,1).
\end{align}
Simple calculation shows the following:
\begin{lem} \label{lem:KL_distance:upper_bound}
Let
$
y=g(\vB^{\tau}\vx)+\epsilon, \epsilon\sim N(0,1)
$
where $\vB \in \mathbb{O}(p,d)$ and $\vx\sim N(0,\vI_{p})$ and let $p_{B,g}(y,\vx)$ be the joint density function of $(y,\vx)$,  then we have
\begin{align}
KL(p_{B,g},p_{B',g}) \leq  |\nabla g|^{2} \|B-B'\|^{2}_{F}.
\end{align}

\begin{comment}
The density function $p_{\vB}(y,\vx)$ of $(y,x)$ is 
\begin{align}
p_{\vB}(y,x)&=p_{\vB}(y|\vx)p(\vx)\\
&=\frac{1}{\sqrt{2\pi}}\exp^{-\frac{1}{2}(y-g(\vB^{\tau}\vx))^{2}}p(\vx)
\end{align}
where $p(\vx)$ is the density function of standard p-dimensional normal distribution.
Let $z=y-g(\vB^{\tau}\vx)$
\begin{align*}
&KL(p_{\vB},p_{\vB'})\\
=&\int \frac{1}{\sqrt{2\pi}}\exp^{-\frac{1}{2}(y-g(\vB^{\tau}\vx))^{2}}p(\vx) \left(\frac{1}{2}\left(y-g(\vB^{'\tau}\vx)\right)^{2}-\frac{1}{2}\left(y-g(\vB^{\tau}\vx)\right)^{2} \right)d\vx dy\\
=&\int \frac{1}{\sqrt{2\pi}}\exp^{-\frac{1}{2}z^{2}}p(\vx)\left(\frac{1}{2}\left(z+g(\vB^{\tau}\vx)-g(\vB^{'\tau}\vx)\right)^{2}-\frac{1}{2}z^{2}  \right)d\vx dz\\
=&\int \frac{1}{\sqrt{2\pi}}\exp^{-\frac{1}{2}z^{2}}p(\vx)\frac{1}{2}\left(g(\vB^{\tau}\vx)-g(\vB^{'\tau}\vx)\right)^{2}d\vx dz\\
\leq & |\nabla g|^{2} \int  p(\vx)\vx^{\tau}\left(\vB-\vB' \right)\left(\vB^{\tau}-\vB^{'\tau} \right)\vx d\vx\\
=&|\nabla g|^{2}Tr\left( \left(\vB-\vB' \right)\left(\vB^{\tau}-\vB^{'\tau} \right)\right)
\end{align*}
\end{comment}
\end{lem}
If $f$ satisfying the Conjecture \ref{conj:derivative},  the Fano Lemma  gives us
\begin{align}
\nonumber &\sup_{\vV \in \widetilde{\Theta}}\bbE\|P_{\widehat{\vV}}-P_{\vV}\|^{2}_{F}\\
\label{eq:fano:argument} \geq& \min \|P_{\vV_{i}}-P_{\vV_{j}}\|^{2}_{F} \left(1-\frac{\max KL(p^{n}_{\vV_{i},f},p^{n}_{\vV_{j},f})+\log(2)}{\log(|\Theta|)}\right)\\
\nonumber \geq & \varepsilon^{2}\left(1-\frac{n\lambda C^{2}_{2}\varepsilon^{2}+4\log2}{\log(|\Theta|)} \right).
\end{align}
Since $\log|\Theta|>Cd(p-d)$, we know that, if
 $\frac{\log(|\Theta|)}{n\lambda}$ is sufficiently small, we have 
\begin{align}
&\sup_{u \in \Theta}\bbE\|P_{\widehat{\vV}}-P_{\vV}\|^{2}_{F}\succ \frac{d(p-d)}{n\lambda}
\end{align}
by choosing $\varepsilon^{2}=\frac{\log(|\Theta|)}{2n\lambda}$.
This gives us the desired lower bound for `Oracle risk'.

\paragraph{Proof of Theorem \ref{thm:rsik:sparse:d:fixed:lambda:fixed}} I. Exact sparsity. With the lower bound of the `Oracle risk',
 we only need to prove the following to obtain the lower bound of the problem with exact sparsity.
\begin{align}
%\label{eq:temp:inline1:nowhere}\inf_{\widehat{\vV}} %\sup_{\mathcal{M}\in \mathfrak{M}_{s,0}\left( p,d,\lambda,%\kappa\right) } \bbE_{\mathcal{M}} \|%\widehat{\vV}\widehat{\vV}^{\tau}-\vV\vV^{\tau}\|^{2}_{F}&\succ %d\wedge \frac{ds}{n\lambda} \\
\label{eq:temp:inline2:nowhere}\inf_{\widehat{\vV}} \sup_{\mathcal{M}\in \mathfrak{M}_{s,0}\left( p,d,\lambda,\kappa\right) } \bbE_{\mathcal{M}} \|\widehat{\vV}\widehat{\vV}^{\tau}-\vV\vV^{\tau}\|^{2}_{F}&\succ d\wedge \frac{s\log\frac{ep}{s}}{n\lambda} .
\end{align}
%\eqref{eq:temp:inline1:nowhere} follows from the lower bound of oracle risk. 
%The \eqref{eq:temp:inline2:nowhere} could be obtain 
It follows from the arguments in \cite{vu2012minimax} and \cite{cai2013sparse}. More precisely,  \cite{vu2012minimax} have constructed a set $\Theta'\subset \bbS^{p-s}$, such that 
\begin{itemize}
\item[1.] $\delta/\sqrt{2}<\|\beta_{1}-\beta_{2}\|_{2}\leq \sqrt{2}\delta$ for all distinct pairs $\beta_{1},\beta_{2} \in \Theta'$,
\item[2.] $\|\beta\|_{0} \leq s$ for all $\beta \in \Theta'$, 
\item[3.] $\log|\Theta'|\geq cs[\log(p-s+1)-\log(s)]$, where $c \geq 0.233$.
\end{itemize}
Now we consider the following family of models
\begin{align*}
y=f(\vV^{\tau}\vx)+\epsilon
\end{align*}
where $\vx \sim N(0,\vI_{p})$, $\epsilon \sim N(0,1)$,
$
\vV=\left(\begin{array}{cc}
\beta & 0_{(p+1-s)\times (s-1)} \\
0_{(s-1)\times 1}& I_{s-1}
\end{array} \right)
$
and $\beta \in \Theta' \subset S^{p-s}$. The similar Fano type argument near \eqref{eq:fano:argument} gives us the \eqref{eq:temp:inline2:nowhere}.

II. Weak $l_{q}$ sparsity.
For the lower bound of problem with  weak $l_{q}$ sparsity, we can simply apply the argument of Theorem 2 in \cite{cai2013sparse}.

\epf

\subsection*{A linear algebraic lemma} In this section, we include a differential geometric argument for the following linear algebraic lemma.
\begin{lem}\label{lem:temp}Let  $\pi : A \mapsto AA^{\tau}$ be the tautological map between the set $\mathbb{O}(p,d) \subset \mathbb{R}^{p\times d}$ of $p\times d$ orthogonal matrices   and  the set  $\mathbb{G}(p,d)\subset \mathbb{R}^{p\times p}$ of all d-dimensional subspaces  in $\mathbb{R}^{p}$. \footnote{$\mathbb{G}(p,d)$ is the so-called Grassmannian. Each point in $\mathbb{G}(p,d)$ can be identified with a projection matrix, thus we can embed  it into $\mathbb{R}^{p^{2}}$.}
There exists an open set $U \subset \mathbb{G}(p,d)$ and positive constant $C$, such that for any $n$ and $u_{1},....,u_{n} \in U$, there exist orthogonal $p\times d$ matrices $A_{1},..,A_{n}$ such that for any $i,j$ one has $A_{j}A_{j}^{\tau}=u_{j}$ and 
\begin{align}
\|A_{i}-A_{j}\|_{F}\leq C\|A_{i}A^{\tau}_{i}-A_{j}A^{\tau}_{j}\|_{F}
\end{align}
\proof Since we did not find an algebraic proof of this lemma, we resort it to the following geometric lemma. the proof of which is rather involved. \epf
\end{lem}

\begin{lem}\label{lem:norm}
The standard metric on Euclidean space induced the Frobenius distance on $\mathbb{O}(p,d)$ by  and the metric on $\mathbb{G}(p,d)$, 
i.e., 
for $A, B \in \mathbb{O}(p,d)$, one has $d(A,B)=\|A-B\|_{F}$ and for $AA^{\tau} , BB^{\tau}  \in \mathbb{G}(p,d)$, one has $D(AA^{\tau},BB^{\tau})=\|AA^{\tau},BB^{\tau}\|_{F}$. 
There exist an open set $U$ and a section $\sigma : U \mapsto  \pi^{-1}(U)$, i.e., $(\pi \circ \sigma =Id_{U})$ in the following commutative diagram, 
\begin{align}
\begin{array}[c]{ccccc}
\mathbb{O}(p,d)                  & \hookleftarrow     & \pi^{-1}(U)               &              & \\
\downarrow \pi                   &                & \downarrow               \pi & \overset{\sigma}{\nwarrow}      & \\
\mathbb{G}(p,d)                  &\hookleftarrow      &U                          &=  & U
\end{array}
\end{align}
satisfying that, there exist positive constants $C_{1}$ and $C_{2}$ such that, for any $A,B \in \sigma(U)$, 
\begin{align}\label{eq:norm}
C_{1}D(\pi(A),\pi(B)) \leq d(A,B)\leq C_{2}D(\pi(A),\pi(B)).
\end{align}

%%%Here end
%Let $I_{1}=\{~i ~|~ i \in S, ~ \bLambda(i,i)>t(1-1/2\nu)^{-1}~\}$,  $I_{2}=\{~i ~|~ i \in S,~ \bLambda(i,i)<t(1+1/2\nu)^{-1}\}$ , $I_{3}=\{~i~ |~ i~ \not \in S ,~\bLambda(i,i)>t\}$ and $I_{4}=\{~i~ |~ i \not \in S,~ \bLambda(i,i)<t~\}$
%
%
%Note that from lemma \ref{lem:main:deviation}, we have for any $i \in T_{t^{2}}$, we have
%\begin{align}
%\bbP(\widehat{\bLambda}_{H}(i,i)>(1-\frac{1}{2\nu})\bLambda(i,i))\leq C_{1}\exp\left(-\frac{n \bLambda(i,i)} {H^{2}}+C_{2}\log(H) \right)
%\end{align} 
%
%\begin{align}
%\bbP(\widehat{\bLambda}_{H}(i,i)<(1+\frac{1}{2\nu})\bLambda(i,i))\leq C_{1}\exp\left(-\frac{n \bLambda(i,i)} {H^{2}}+C_{2}\log(H) \right)
%\end{align} 
%
%\begin{align}
%\bbP(\widehat{\bLambda}_{H}(i,i)>
%\end{align}

\proof
\footnote{To avoid reproducing standard content in textbook, in this proof, we assume the reader has some familiarities  with differential geometry .} 
Note that for any submanifold $M \subset \mathbb{R}^{n}$, there are two distances: the induced distance $d_{M}^{i}$ which induced from the Euclidean distance  and the geodesic distance $d_{M}^{g}$ where $M$ with the induced Riemannian metric.
\begin{lem} We embed $\mathbb{O}(p,d)$ into the set of all $p\times d$ matrices
and $\mathbb{G}(p,d)$ into the set of all $p\times p$ matrices.  There exist constants $C_{1}$, $C_{2}$ such that, for any point $A \in \mathbb{O}(p,d)$,  there exists an open neighbourhood $U$ of $A$ such that  for any $u_{1}, u_{2} \in U$, one has
\begin{align}\label{eq:different:norm1}
C_{1}d_{\mathbb{O}(p,d)}^{i}(u_{1},u_{2}) \leq d_{\mathbb{O}(p,d)}^{g}(u_{1},u_{2})\leq C_{2}d_{\mathbb{O}(p,d)}^{i}(u_{1},u_{2}).
\end{align}
Similarly, there exist constants $C_{1}$, $C_{2}$ such that, for any point $S_{A} \in \mathbb{G}(p,d)$,  there exists an open neighbourhood $U$ of $S_{A}$ such that for any $u_{1}, u_{2} \in U$, one has
\begin{align}\label{eq:different:norm2}
C_{1}d_{\mathbb{G}(p,d)}^{i}(u_{1},u_{2}) \leq d_{\mathbb{G}(p,d)}^{g}(u_{1},u_{2})\leq C_{2}d_{\mathbb{G}(p,d)}^{i}(u_{1},u_{2}).
\end{align}
\proof We only prove the inequality \eqref{eq:different:norm1}. The first part is trivial. Thus we only need to prove the second inequality. Let $J=(I_{d},0_{d\times(p-d)})^{\tau}$. Since $\mathbb{O}(p,d)$ is homogeneous manifold, we only need to prove that there exists a neighbourhood $U$ of $J$ such that for any $u \in U$, the \eqref{eq:different:norm1} holds for $J$ and $u$. Note that any orthogonal matrix $u$ near $J$ could be written as
\begin{align}
u=\left(\begin{array}{c}
I_{d}\\
B
\end{array}\right)(1+B^{\tau}B)^{-1/2}T
\end{align}
where $B$ is some $(p-d)\times d$ matrix and $T$ is an $d\times d$ orthogonal matrix. 

Let us consider the curve 
\begin{align}
u(t)=\left(\begin{array}{c}
I_{d}\\
tB
\end{array}\right)(1+t^{2}B^{\tau}B)^{-1/2}T
\end{align}
inside $\mathbb{O}(p,d)$ and denote its length by $s$. Since $d^{g}_{\mathbb{O}(p,d)}(J,u) \leq s$, we only need to prove that there exists constant $C$ such that $s \leq C\|J-u\|_{F}$. 
Note that
\begin{align}
\frac{du}{dt}=\left(\begin{array}{c}
0\\
B
\end{array}\right)(1+t^{2}B^{\tau}B)^{-1/2}T
-
t\left(\begin{array}{c}
I_{d}\\
tB
\end{array}\right)(1+t^{2}B^{\tau}B)^{-3/2}B^{\tau}BT
\end{align}
Let $B^{\tau}B=VE^{2}V^{\tau}$ where $E^{2}$ is diagonal matrix with entries $e^{2}_{1},...,e^{2}_{d}$. Then one has
\begin{align}
\|\frac{du}{dt}\|^{2}
=&tr (E^{2}(1+t^{2}E^{2})^{-2})
\end{align}

%\begin{align}
%\|\frac{du}{dt}\|^{2}=&tr\left(\frac{du}{dt}(\frac{du}{dt})^{\tau}\right)\\
%=&tr \left( (1+t^{2}B^{\tau}B)^{-1}B^{\tau}B\right)
%+t^{2}tr\left(B^{\tau}B(1+t^{2}B^{\tau}B)^{-2}B^{\tau}B \right)\\
%&-2t^{2}tr \left( (1+t^{2}B^{\tau}B)^{-1/2}
%B^{\tau}B(1+t^{2}B^{\tau}B)^{-3/2}B^{\tau}B \right)\\
%=&tr\left(D^{2}(1+t^{2}D^{2})^{-1}+t^{2}D^{4}(1+t^{2}D^{2})^{-2}-2t^{2}D^{4}(1+t^{2}D^{2})^{-2} \right)\\
%=&tr (D^{2}(1+t^{2}D^{2})^{-2})
%\end{align}
For any orthogonal matrix T and semi-positive definite matrix A, one has $tr(TA)\leq tr(A)$. Thus we have
\begin{align*}
\|J-u(t)\|^{2}_{F}=&2d-2tr\left((1+t^{2}B^{\tau}B)^{-1/2}T\right)\geq 2d-2tr\left((1+t^{2}B^{\tau}B)^{-1/2}\right)\\
=&2d-2tr\left((1+t^{2}E^{2})^{-1/2}\right).
\end{align*}

Since
\begin{align}
s^{2}=\left(\int_{0}^{t}\|\frac{du}{da}\|da\right)^{2} \leq t \int_{0}^{t}\|\frac{du}{da}\|^{2}da\leq t \int_{0}^{t}tr (E^{2}(1+t^{2}E^{2})^{-2})da,
\end{align}
we only need to prove that there exists positive constant $C$,  such that for any $1\leq j \leq d$, one has
\begin{align}
t\int^{t}_{0}\frac{e_{j}^{2}}{(1+a^{2}e_{j}^{2})^{2}}da \leq C\left( 1-\frac{1}{\sqrt{1+t^{2}e_{j}^{2}}} \right)
\end{align}
which can be verified directly when both $e_{j}~'s$ and $t$ are sufficiently small.
\epf
\end{lem}

Note that $\mathbb{O}(p,d)$ is a principal bundle over $\mathbb{G}(p,d)$ with structure group $\mathbb{O}(d,d)$. 
Thus  for any point $A \in \mathbb{O}(p,d)$, we have a natural decomposition of tangent space at $A$:
\begin{align}\label{eq:orthogonal:decomposition}
T_{A}\mathbb{O}(p,d)=V_{A}\oplus H_{A}
\end{align}
where 
$
V_{A}=\{~AX~ | X \mbox{ is } d\times d \mbox{ anti-symmetric matrix.}~\}$
consists of vector tangent to the fibre and
$
H_{A}=\{~B~ | B \mbox{ is } p\times d \mbox { matrix such that } B^{\tau}A=0.~\}
$
consists of vector perpendicular to $V_{A}$. For any non-zero $\alpha \in T_{A} \mathbb{O}(p,d)$, let $\alpha=\alpha_{V}\oplus\alpha_{H}$ be the decomposition with respect to \eqref{eq:orthogonal:decomposition}. We introduce 
\begin{align}
\mu(\alpha)\triangleq \frac{\|\alpha_{H}\|}{\|\alpha\|}.
\end{align}
For a subspace $V \subset T_{A}\mathbb{O}(p,d)$, we define 
$\mu(V)=\inf_{\alpha\in V\backslash 0}\mu(\alpha)$.

\begin{lem}
For any $A \in \mathbb{O}(p,d)$, there exists a neighbourhood $U$ of $\pi(A) \in \mathbb{G}(p,d)$ such that there exists a smooth map $\sigma : U \mapsto \pi^{-1}(U)$ satisfying $\sigma(\pi(A))=A$,
$d\sigma|_{\pi(A)}\left(T_{\pi(A)}\mathbb{G}(p,d)\right)=H_{A}$
and 
\begin{align}\label{inlince:neighbour:nondegeneracy}
\forall u \in U, \mu\left(d\sigma|_{u}T_{u}\mathbb{G}(p,d)\right)\geq \frac{1}{2}.
\end{align}
As a direct corollary, we know that there exists two positive constant $C_{1}$, $C_{2}$ such that for any $u_{1}, u_{2} \in U$, one has
\begin{align}\label{inlince:norm:riemanian}
C_{1}d_{\mathbb{G}(p,d)}^{g}(u_{1},u_{2}) \leq d_{\mathbb{O}(p,d)}^{g}(\sigma(u_{1}), \sigma(u_{2}))\leq C_{2}d_{\mathbb{G}(p,d)}^{g}(u_{1},u_{2}).
\end{align}

\proof The existence of U and the inequality \eqref{inlince:neighbour:nondegeneracy} follows from a typical continuity argument.
For the second inequality in \eqref{inlince:norm:riemanian}, let $u(t)$ be a geodesic connect $u_{1}$ and $u_{2}$. When $U$ is sufficient small, $u(t)$ is unique and and 
\begin{align}
d^{g}_{\mathbb{G}(p,d)}(u_{1},u_{2})=\int^{1}_{0}\|\frac{du(t)}{dt}\|_{\mathbb{G}}dt.
\end{align}
Let $L$ be the length of the curve $\sigma(u(t))$, then
\begin{align}
\nonumber L=&\int^{1}_{0}\|\frac{d\sigma(u(t))}{dt}\|_{\mathbb{O}}dt=\int^{1}_{0}\|d\sigma(\frac{du(t)}{dt})\|_\mathbb{O}dt\\
\label{continuity:argument}\leq&2\int_{0}^{1}\|\left(d\sigma\left(\frac{du(t)}{dt}\right)\right)_{H}\|_{\mathbb{O}}dt\\
\label{norm:preserving}= & 2\sqrt{2}\int_{0}^{1} \|\frac{du(t)}{dt}\|_{\mathbb{G}}dt\\
\nonumber= & 2\sqrt{2}d^{g}_{\mathbb{G}(p,d)}(u_{1},u_{2})
\end{align}
where inequality \eqref{continuity:argument} follows from the inequality \eqref{inlince:neighbour:nondegeneracy} and
 equation \eqref{inlince:norm:riemanian} follows from the fact that for any $B \in T_{B}\mathbb{O}(p,d)$,( i.e., $B^{\tau}A+A^{\tau}B=0$), one has 
$
d\pi(B)=AB^{\tau}+BA^{\tau} \in T_{AA^{\tau}}\mathbb{G}(p,d)$
and if $B \in H_{A}$,(i.e., $B^{\tau}A=0$), one has
$
 tr d\pi(B)d\pi(B)^{\tau}=2trBB^{\tau}.
$
In particular, we know that
\begin{align}
d^{g}_{\mathbb{O}(p,d)}(\sigma(u_{1}),\sigma(u_{2}))\leq L \leq 2\sqrt{2} d^{g}_{\mathbb{G}(p,d)}(u_{1},u_{2})
\end{align}
The first inequality in \eqref{inlince:norm:riemanian} can be proved similar and thus omitted.
\epf
\end{lem}

The Lemma \ref{lem:norm} is a direct corollary of the above two Lemmas.\epf
\end{lem}

\subsection*{Assisting Lemmas}
The following lemmas are borrowed from  \cite{vershynin2010introduction} and \cite{cai2013sparse}.
\begin{lem} \label{random:nonasymptotic}
Let $\vE_{p\times H}$ be a $p\times H$ matrix, whose entries are independent standard normal random variables. Then for every $ t\geq 0$, with probability at least $1-2\exp(-t^{2}/2)$, one has :
\[ \lambda^{+}_{sing,min}(\vE_{p\times H}) \geq \sqrt{p}-\sqrt{H}-t\],
and
\[ \lambda_{sing,max}(\vE_{p\times H}) \leq \sqrt{p}+\sqrt{H}+t.\]
\end{lem} 

\begin{cor}\label{cor:final:done}
One has
\[
\frac{1}{2}\left(\sqrt{p}-\sqrt{H}\right) \leq \lambda^{-}_{sing,min}\left(\vE_{p\times H}\right) \leq \lambda_{sing,max}\left(\vE_{p\times H}\right)\leq \frac{3}{2}\left(\sqrt{p}+\sqrt{H}\right).
\]
with probability converging to one, as $n \rightarrow \infty$.
\end{cor}

\begin{lem} \label{lem:sin_theta} 
$($Sin-Theta Theorem. $)$ Let $\vA$ and $\vA+\vE$ be symmetric matrices satisfying
\begin{equation*}
\vA=[\vF_{0},\vF_{1}]\left[ \begin{array}{cc}
\vA_{0} & 0\\
0 & \vA_{1}
\end{array} \right]
\left[ \begin{array}{c}
\vF^{\tau}_{0}\\
\vF^{\tau}_{1}
\end{array} \right]
\quad
\vA+\vE=[\vG_{0},\vG_{1}]\left[ \begin{array}{cc}
\bLambda_{0} & 0\\
0 & \bLambda_{1}
\end{array} \right]
\left[ \begin{array}{c}
\vG^{\tau}_{0}\\
\vG^{\tau}_{1}
\end{array} \right]
\end{equation*}
where $[\vF_{0},\vF_{1}]$ and $[\vG_{0},\vG_{1}]$ are orthogonal matrices. If the eigenvalues of $\vA_{0}$ are contained in an interval (a,b) , and the eigenvalues of $\bLambda_{1}$ are excluded from the interval $(a-\delta,b+\delta)$ for some $\delta>0$, then
\[
\|\vF_{0}\vF_{0}^{\tau}-\vG_{0}\vG_{0}^{\tau}\|\leq \frac{\min(\|\vF_{1}^{\tau}\vE\vG_{0}\|, \|\vF_{0}^{\tau}\vE\vG_{1}\|)}{\delta},
\]
and 
\[
\frac{1}{\sqrt{2}}\|\vF_{0}\vF_{0}^{\tau}-\vG_{0}\vG_{0}^{\tau}\|_{F} \leq \frac{\min(\|\vF_{1}^{\tau}\vE\vG_{0}\|_{F}, \|\vF_{0}^{\tau}\vE\vG_{1}\|_{F})}{\delta}.
\]
\end{lem}

\begin{lem}\label{lem:Cai:lem4}
Let $\vK \in \bbR^{p\times p}$ be symmetric such that $Tr(\vK)=0$ and $\|\vK\|_{F}\leq1$. Let $\vZ$ be an $H \times p$  matrix consisting of independent standard normal entries. Then for any $t>0$, one has
\begin{align}
\bbP\left(\Big|\Big\langle \vZ^{\tau}\vZ, \vK \Big\rangle\Big
|\geq 2\sqrt{H}t +2t^{2} \right) \leq 2\exp\left(-t^{2} \right).
\end{align}
\end{lem}
We remind that this lemma is a trivial modification of  Lemma 4 in \cite{cai2013sparse}, where they assumed $\|\vK\|_{F}=1$. 

\begin{lem}\label{lem:Cai:lem5}
Let $X_{1},...,X_{N}$ be i.i.d such that
\begin{align}
\bbP(|X_{1}|\geq at+bt^{2})\leq c\exp\left(-t^{2} \right)
\end{align}
where $a,b,c>0$. Then
\begin{align}
\bbE\max|X_{i}|^{2} \leq (2a^{2}+8b^{2})\log(ecN)+2b^{2}\log^{2}(cN).
\end{align}
\end{lem}

\end{appendix}

%\bibliographystyle{plainnat}
%\bibliography{sir}
\end{document}